\documentclass[11pt]{article}
\usepackage{graphicx}
\usepackage[utf8]{inputenc}
\usepackage[full]{textcomp}
\usepackage[osf]{newtxtext}
\usepackage{comment}

\usepackage{amssymb}
\usepackage{mathtools}
\usepackage{hyperref}
\usepackage{breakurl}
\usepackage{mhenvs}
\usepackage{mhequ}
\usepackage{mhsymb}
\usepackage{tikz}
\usepackage{mathrsfs}
\usepackage{halloweenmath}
\usepackage{cprotect}

\usepackage{microtype}
\usepackage{wasysym}
\usepackage{centernot}
\usepackage{enumitem}

\usepackage{geometry}

\usetikzlibrary{shapes.misc}
\usetikzlibrary{shapes.symbols}
\usetikzlibrary{shapes.geometric}
\usetikzlibrary{decorations}
\usetikzlibrary{decorations.markings}
\usetikzlibrary{snakes} 

\usetikzlibrary{calc}
\usetikzlibrary{external}

\let\oldskull\skull
\def\skull{\mathord{\oldskull}}
	
\DeclareMathAlphabet{\mathbbm}{U}{bbm}{m}{n}

\DeclareFontFamily{U}{BOONDOX-calo}{\skewchar\font=45 }
\DeclareFontShape{U}{BOONDOX-calo}{m}{n}{
  <-> s*[1.05] BOONDOX-r-calo}{}
\DeclareFontShape{U}{BOONDOX-calo}{b}{n}{
  <-> s*[1.05] BOONDOX-b-calo}{}
\DeclareMathAlphabet{\mcb}{U}{BOONDOX-calo}{m}{n}
\SetMathAlphabet{\mcb}{bold}{U}{BOONDOX-calo}{b}{n}

\setlist{noitemsep,topsep=4pt}

\makeatletter
\def\DeclareSymbol#1#2#3{%
	\expandafter\gdef\csname MH@symb@#1\endcsname{\tikzsetnextfilename{symbol#1}%
	\tikz[baseline=#2,scale=0.15,draw=symbols,line join=round,line cap=round]{#3}}%
	\expandafter\gdef\csname MH@symb@#1s\endcsname{\scalebox{0.75}{\tikzsetnextfilename{symbol#1}%
	\tikz[baseline=#2,scale=0.15,draw=symbols,line join=round,line cap=round]{#3}}}%
	\expandafter\gdef\csname MH@symb@#1ss\endcsname{\scalebox{0.65}{\tikzsetnextfilename{symbol#1}%
	\tikz[baseline=#2,scale=0.15,draw=symbols,line join=round,line cap=round]{#3}}}%
	}
\def\<#1>{\ifthenelse{\boolean{mmode}}{\mathchoice{\csname MH@symb@#1\endcsname}{\csname MH@symb@#1\endcsname}{\csname MH@symb@#1s\endcsname}{\csname MH@symb@#1ss\endcsname}}{\csname MH@symb@#1\endcsname}}
\makeatother

\DeclareSymbol{1}{0}{\draw[white] (-.5,0) -- (.5,0); \draw (0,0)  -- (0,1.5) node[sdot] {};}
\DeclareSymbol{11}{0}{\draw (-0.5,1.2) node[sdot] {} -- (0,0) -- (0.5,1.2) node[sdot] {};}
\DeclareSymbol{111}{0}{\draw (0,0) -- (0,1.2) node[sdot] {}; \draw (-.7,1) node[sdot] {} -- (0,0) -- (.7,1) node[sdot] {};}
\DeclareSymbol{131}{-3}{\draw (0,0) -- (0,-1) -- (1,0) node[sdot] {}; \draw (0,0) -- (0,-1) -- (-1,0) node[sdot] {}; \draw (0,0) -- (0,1.2) node[sdot] {}; \draw (-.7,1) node[sdot] {} -- (0,0) -- (.7,1) node[sdot] {};}
\DeclareSymbol{11}{0}{\draw (-0.5,1.2) node[sdot] {} -- (0,0) -- (0.5,1.2) node[sdot] {};}
\DeclareSymbol{12}{-3}{\draw (-.8,1) node[sdot] {} -- (0,0) -- (0,-1); \draw (1,0) node[sdot] {} -- (0,-1) -- (-1,0) node[sdot] {};}
\DeclareSymbol{30}{-3}{\draw (0,0) -- (0,-1); \draw (0,0) -- (0,1.2) node[sdot] {}; \draw (-.7,1) node[sdot] {} -- (0,0) -- (.7,1) node[sdot] {};}
\DeclareSymbol{10}{-3}{\draw (-.8,1) node[sdot] {} -- (0,0) -- (0,-1);}
\DeclareSymbol{22}{-3}{\draw (0,0.3) -- (0,-1) -- (1,0) node[sdot] {}; \draw (0,0.3) -- (0,-1) -- (-1,0) node[sdot] {};\draw (-.7,1) node[sdot] {} -- (0,0.3) -- (.7,1) node[sdot] {};}
\DeclareSymbol{211}{0}{\draw (-0.5,2.4) node[sdot] {} -- (-1,1.6);\draw (-1.5,2.4) node[sdot] {} -- (-0.5,0.8); \draw (0,1.6) node[sdot] {} -- (-0.5,0.8) -- (0,0) -- (0.5,0.8) node[sdot] {};}

\DeclareSymbol{Xi22}{0.5}{\draw (0,0) node[xi] {} -- (-1,1) node[not] {} -- (0,2) node[xi] {};}

\DeclareSymbol{Xi2}{-2}{\draw (0,-0.25) node[xi] {} -- (-1,1) node[xi] {};}
\DeclareSymbol{Xi3}{0}{\draw (0,0) node[xi] {} -- (-1,1) node[xi] {} -- (0,2) node[xi] {};}
\DeclareSymbol{Xi4}{2}{\draw (0,0) node[not] {} -- (-1,1) node[xi] {} -- (0,2) node[xi] {} -- (-1,3) node[xi] {};}
\DeclareSymbol{Xi2X}{-2}{\draw (0,-0.25) node[xi] {} -- (-1,1) node[xix] {};}
\DeclareSymbol{XXi2}{-2}{\draw (0,-0.25) node[xix] {} -- (-1,1) node[xi] {};}

\DeclareSymbol{IXi^2asym}{-1}{\draw (-1,1) node[xired] {} -- (0,-.4);
\draw (1,1) node[xigreen] {} -- (0,-.4);}

\DeclareSymbol{IXi^2asym2}{-1}{\draw (-1,1) node[xigreen] {} -- (0,-.4);
\draw (1,1) node[xired] {} -- (0,-.4);}

\DeclareSymbol{IXi^2asym3}{-1}{\draw (-1,1) node[xired] {} -- (0,-.4);
\draw (1,1) node[xigreen] {} -- (0,-.4);}

\DeclareSymbol{IXi^2}{-1}{\draw (-1,1) node[xi] {} -- (0,-.4);
\draw (1,1) node[xi] {} -- (0,-.4);}

\DeclareSymbol{rIXi^2}{-1}{\draw (-1,1) node[xired] {} -- (0,-.4);
\draw (1,1) node[xired] {} -- (0,-.4);}

\DeclareSymbol{IXi^2green}{-1}{\draw (-1,1) node[xigreen] {} -- (0,-.4);
\draw (1,1) node[xigreen] {} -- (0,-.4);}

\DeclareSymbol{IXi^2_typed}{-1}{\draw (-1,1) node[xiorange] {} -- (0,-.4);
\draw (1,1) node[xigreen] {} -- (0,-.4);}

\DeclareSymbol{IXi^2green*}{-1}{\draw (-1,1) node[xigreen1] {} -- (0,-.4);
\draw (1,1) node[xigreen1] {} -- (0,-.4);}

\DeclareSymbol{IXiI'Xi_notriangle}{-1}{\draw (-1,1) node[xi] {} -- (0,0);
\draw (1,1)[kernels2] node[xi] {} -- (0,0) {};}

\DeclareSymbol{IXiI'Xi}{-1}{\draw (-1,1) node[xi] {} -- (0,0);
\draw (1,1)[kernels2] node[xi] {} -- (0,0) node[not] {};}

\DeclareSymbol{IXiI'Xi_typed1}{-1}{\draw (-1,1) node[xigreen] {} -- (0,0);
\draw (1,1)[kernels2,greennode] node[xired] {} -- (0,0) node[not] {};}

\DeclareSymbol{IXiI'Xi_typed2}{-1}{\draw (-1,1) node[xigreen] {} -- (0,0);
\draw (1,1)[kernels2,rednode] node[xigreen] {} -- (0,0) node[not] {};}

\DeclareSymbol{IXiI'Xi_typed2_notriangle}{-1}{\draw (-1,1) node[xigreen] {} -- (0,0);
\draw (1,1)[kernels2,rednode] node[xigreen] {} -- (0,0) {};}

\DeclareSymbol{IXi^3}{-1}{\draw (-1.3,1) node[xi] {} -- (0,0);
\draw (1.3,1) node[xi] {} -- (0,0);
\draw (0,1.5) node[xi] {} -- (0,0) node[not] {};}

\DeclareSymbol{IXi^3_notriangle}{0}{\draw (-1.3,1) node[xi] {} -- (0,0);
\draw (1.3,1) node[xi] {} -- (0,0);
\draw (0,1.5) node[xi] {} -- (0,0) {};}

\DeclareSymbol{IXi^3_typed}{-1}{\draw (-1.3,1) node[xigreen] {} -- (0,0);
\draw (1.3,1) node[xired] {} -- (0,0);
\draw (0,1) node[xigreen] {} -- (0,0) node[not] {};}

\DeclareSymbol{IXi^3-typed112}{-1}{\draw (-1.3,1) node[xi] {} -- (0,0);
\draw (1.3,1) node[xigreen] {} -- (0,0);
\draw (0,1) node[xi] {} -- (0,0) node[not] {};}

\DeclareSymbol{I'Xi}{-1}{
\draw[kernels2] (1,1) node[xi] {} -- (0,0) node[not] {};}

\DeclareSymbol{I'Xi_notriangle}{-1}{
\draw[kernels2] (1,1) node[xi] {} -- (0,0) {};}

\DeclareSymbol{I'Xigreen-red}{-1}{
\draw[kernels2,rednode] (0,1) node[xi,fill=greennode] {} -- (0,-.4);}

\DeclareSymbol{IXi}{-1}{
\draw (0,1) node[xi] {} -- (0,-.4);}

\DeclareSymbol{IXigreen}{-1}{
\draw (0,1) node[xigreen] {} -- (0,-.4);}

\DeclareSymbol{IXired}{-1}{
\draw (0,1) node[xired] {} -- (0,-.4);}

\DeclareSymbol{green}{-1}{
\draw[greennode,very thick,line cap=round,scale=1.3] (0.05,0.4) -- (0.2,0.55) -- (0,0) -- (0.6,0.6) -- (0.4,0.05) -- (0.55,0.2);}
\DeclareSymbol{red}{-1}{
\draw[rednode,very thick,line cap=round,scale=1.3] (0.05,0.4) -- (0.2,0.55) -- (0,0) -- (0.6,0.6) -- (0.4,0.05) -- (0.55,0.2);}
\DeclareSymbol{lblue}{-1}{
\draw[lbluenode,very thick,line cap=round,scale=1.3] (0.05,0.4) -- (0.2,0.55) -- (0,0) -- (0.6,0.6) -- (0.4,0.05) -- (0.55,0.2);}
\DeclareSymbol{dblue}{-1}{
\draw[dbluenode,very thick,line cap=round,scale=1.3] (0.05,0.4) -- (0.2,0.55) -- (0,0) -- (0.6,0.6) -- (0.4,0.05) -- (0.55,0.2);}
\DeclareSymbol{orange}{-1}{
\draw[orangenode,very thick,line cap=round,scale=1.3] (0.05,0.4) -- (0.2,0.55) -- (0,0) -- (0.6,0.6) -- (0.4,0.05) -- (0.55,0.2);}

\DeclareSymbol{I'XXi}{-1}{
\draw[kernels2] (1,1) node[xi] {} -- (0,0) node[not] {};
\draw (1,1) node[xix] {};}

\DeclareSymbol{I'XXi_notriangle}{-1}{
\draw[kernels2] (1,1) node[xi] {} -- (0,0) {};
\draw (1,1) node[xix] {};}

\DeclareSymbol{IXiI'[IXiI'Xi]}{-1}{\draw[kernels2] (2,2) node[xi] {} -- (1,1);
\draw (0,2) node[xi] {} -- (1,1) node[not] {};
\draw[kernels2] (1,1) -- (0,0) node[not] {};
\draw (-1,1) node[xi] {} -- (0,0);}

\DeclareSymbol{IXiI'[IXiI'Xi]_notriangle}{1}{\draw[kernels2] (2,2) node[xi] {} -- (1,1);
\draw (0,2) node[xi] {} -- (1,1) {};
\draw[kernels2] (1,1) -- (0,0) {};
\draw (-1,1) node[xi] {} -- (0,0);}

\DeclareSymbol{I'[IXiI'Xi]}{-1}{\draw[kernels2] (2,2) node[xi] {} -- (1,1);
\draw (0,2) node[xi] {} -- (1,1) node[not] {};
\draw[kernels2] (1,1) -- (0,0) node[not] {};}

\DeclareSymbol{I'[IXiI'Xi]_notriangle}{1}{\draw[kernels2] (2,2) node[xi] {} -- (1,1);
\draw (0,2) node[xi] {} -- (1,1) {};
\draw[kernels2] (1,1) -- (0,0) {};}

\DeclareSymbol{I[IXi^2]IXi-typed112}{-1}{\draw (-2,2) node[xi] {} -- (-1,1);
\draw (0,2) node[xi] {} -- (-1,1) node[not] {};
\draw (-1,1) -- (0,0) node[not] {};
\draw (1,1) node[xigreen] {} -- (0,0);}

\DeclareSymbol{I[IXi^2]IXi-typed211}{-1}{\draw (-2,2) node[xigreen] {} -- (-1,1);
\draw (0,2) node[xi] {} -- (-1,1) node[not] {};
\draw (-1,1) -- (0,0) node[not] {};
\draw (1,1) node[xi] {} -- (0,0);}

\DeclareSymbol{I[IXiI'Xi]I'Xi}{-1}{\draw (-2,2) node[xi] {} -- (-1,1);
\draw[kernels2] (0,2) node[xi] {} -- (-1,1) node[not] {};
\draw (-1,1) -- (0,0) node[not] {};
\draw[kernels2] (1,1) node[xi] {} -- (0,0);}

\DeclareSymbol{I[IXiI'Xi]I'Xi_notriangle}{1}{\draw (-2,2) node[xi] {} -- (-1,1);
\draw[kernels2] (0,2) node[xi] {} -- (-1,1) {};
\draw (-1,1) -- (0,0) {};
\draw[kernels2] (1,1) node[xi] {} -- (0,0);}

\DeclareSymbol{I[I[Xi]]Xi_notriangle}{1}{
\draw (0,2) node[xi] {} -- (-1,1) {};
\draw[snake=snake , segment length=1pt, segment amplitude=1pt] (-1,1) -- (0,0) node[xi] {};}

\DeclareSymbol{I[I[Xi]]Xi_typed}{1}{
\draw (0,2) node[xigreen] {} -- (-1,1) {};
\draw (-1,1) -- (0,0) node[xigreen] {};}

\DeclareSymbol{I[I[Xi]]Xi_typed*}{1}{
\draw[snake=snake , segment length=1pt, segment amplitude=1pt] (0,2) node[xigreen1] {} -- (-1,1) {};
\draw (-1,1) -- (0,0) node[xigreen1] {};}

\DeclareSymbol{I[I'Xi]I'Xi}{1}{\draw[kernels2] (0,2) node[xi] {} -- (-1,1);
\draw (-1,1) node[fill=rednode,not] {} -- (0,0);
\draw[kernels2] (1,1) node[xi] {} -- (0,0);}

\DeclareSymbol{I[I'Xi]I'Xi_notriangle}{1}{\draw[kernels2] (0,2) node[xi] {} -- (-1,1);
\draw (-1,1) {} -- (0,0);
\draw[kernels2] (1,1) node[xi] {} -- (0,0);}

\DeclareSymbol{I[I'Xi]I'Xi_typed}{1}{\draw[kernels2,rednode] (0,2) node[xigreen] {} -- (-1,1);
\draw (-1,1) node[notorange] {} -- (0,0);
\draw[kernels2,rednode] (1,1) node[xigreen] {} -- (0,0);}

\DeclareSymbol{I[I'Xi]I'Xi_typed*}{1}{\draw[kernels2,rednode] (0,2) node[xigreen1] {} -- (-1,1);
\draw (-1,1) node[notorange] {} -- (0,0);
\draw[kernels2,rednode] (1,1) node[xigreen1] {} -- (0,0);}

\DeclareSymbol{I[I'Xi]I'Xi_typed-h}{1}{\draw[kernels2,greennode] (0,2) node[xigreen] {} -- (-1,1);
\draw[snake=snake , segment length=3pt, segment amplitude=1pt] (-1,1) node[notgreen] {} -- (0,0);
\draw[kernels2,greennode] (1,1) node[xigreen] {} -- (0,0);}

\DeclareSymbol{I[I'Xi]I'Xi_typed-h*}{1}{\draw[kernels2,greennode] (0,2) node[xigreen1] {} -- (-1,1);
\draw[snake=snake , segment length=3pt, segment amplitude=1pt] (-1,1) node[notgreen] {} -- (0,0);
\draw[kernels2,greennode] (1,1) node[xigreen1] {} -- (0,0);}

\DeclareSymbol{I[I'Xi]I'Xib}{1}{\draw[kernels2] (0,2) node[xi] {} -- (1,1);
\draw (1,1) node[not] {} -- (0,0);
\draw[kernels2] (-1,1) node[xi] {} -- (0,0);}

\DeclareSymbol{I[I'Xi]I'Xib_notriangle}{1}{\draw[kernels2] (0,2) node[xi] {} -- (1,1);
\draw (1,1) {} -- (0,0);
\draw[kernels2] (-1,1) node[xi] {} -- (0,0);}

\DeclareSymbol{IXiI'[I'Xi]}{1}{\draw[kernels2] (0,2) node[xi] {} -- (1,1);
\draw[kernels2] (1,1) node[not] {} -- (0,0);
\draw (-1,1) node[xi] {} -- (0,0);}

\DeclareSymbol{IXiI'[I'Xi]_notriangle}{1}{\draw[kernels2] (0,2) node[xi] {} -- (1,1);
\draw[kernels2] (1,1)  {} -- (0,0);
\draw (-1,1) node[xi] {} -- (0,0);}

\DeclareSymbol{IXiI'[I'Xi]_typed}{1}{\draw[kernels2,rednode] (0,2) node[xigreen] {} -- (1,1);
\draw[kernels2,rednode] (1,1) node[notorange] {} -- (0,0);
\draw (-1,1) node[xigreen] {} -- (0,0);}

\DeclareSymbol{IXiI'[I'Xi]_typed-h}{1}{\draw[kernels2,greennode] (0,2) node[xigreen] {} -- (1,1);
\draw[kernels2,greennode,snake=snake , segment length=3pt, segment amplitude=1pt] (1,1) node[notgreen] {} -- (0,0);
\draw (-1,1) node[xigreen] {} -- (0,0);}

\DeclareSymbol{IXiI'[I'Xi]_typed*}{1}{\draw[kernels2,rednode] (0,2) node[xigreen1] {} -- (1,1);
\draw[kernels2,rednode] (1,1) node[notorange] {} -- (0,0);
\draw (-1,1) node[xigreen1] {} -- (0,0);}

\DeclareSymbol{IXiI'[I'Xi]_typed-h*}{1}{\draw[kernels2,greennode] (0,2) node[xigreen1] {} -- (1,1);
\draw[kernels2,greennode,snake=snake , segment length=3pt, segment amplitude=1pt] (1,1) node[notgreen] {} -- (0,0);
\draw (-1,1) node[xigreen1] {} -- (0,0);}

\DeclareSymbol{I[I'Xi]}{-1}{\draw[kernels2] (0,2) node[xi] {} -- (1,1);
\draw (1,1) node[not] {} -- (0,0) node[not] {};}

\DeclareSymbol{I[I'Xi]_typed}{1}{\draw[kernels2,rednode] (0,2) node[xigreen] {} -- (-1,1);
\draw (-1,1) node[notorange] {} -- (0,0);}

\DeclareSymbol{I'[I'Xi]}{-1}{\draw[kernels2] (0,2) node[xi] {} -- (1,1);
\draw[kernels2] (1,1) node[not] {} -- (0,0) node[not] {};}

\DeclareSymbol{I'[I'Xi]_notriangle}{-1}{\draw[kernels2] (0,2) node[xi] {} -- (1,1);
\draw[kernels2] (1,1) {} -- (0,0) {};}

\DeclareSymbol{I'[I'Xi]_typed}{1}{\draw[kernels2,rednode] (0,2) node[xigreen] {} -- (1,1);
\draw[kernels2,rednode] (1,1) node[notorange] {} -- (0,0)  {};}

\DeclareSymbol{IXiI'XXi}{-1}{\draw (-1,1) node[xi] {} -- (0,0);
\draw[kernels2] (1,1) node[xi] {} -- (0,0) node[not] {};
\draw (1,1) node[xix] {};}

\DeclareSymbol{IXiI'XXi_notriangle}{-1}{\draw (-1,1) node[xi] {} -- (0,0);
\draw[kernels2] (1,1) node[xi] {} -- (0,0) {};
\draw (1,1) node[xix] {};}

\DeclareSymbol{IXiI'XXi_typed}{-1}{\draw (-1,1) node[xigreen] {} -- (0,-0.3);
\draw[kernels2,rednode] (1,1)  -- (0,-0.3) {};
\draw (1,1) node[xix-green-red] {};
}

\DeclareSymbol{IXiI'XXi_typed*}{-1}{\draw (-1,1) node[xigreen1] {} -- (0,-0.3);
\draw[kernels2,rednode] (1,1)  -- (0,-0.3) {};
\draw (1,1) node[xix-green-red1] {};
}

\DeclareSymbol{IXXiI'Xi}{-1}{\draw (-1,1) node[xix] {} -- (0,0);
\draw[kernels2] (1,1) node[xi] {} -- (0,0) node[not] {};}

\DeclareSymbol{IXXiI'Xi_notriangle}{-1}{\draw (-1,1) node[xix] {} -- (0,0);
\draw[kernels2] (1,1) node[xi] {} -- (0,0) {};}

\DeclareSymbol{IXXiI'Xi_typed}{-1}{\draw (-1,1) node[xix-green-red] {} -- (0,-0.3);
\draw[kernels2,rednode] (1,1) node[xigreen] {} -- (0,-0.3) {};}

\DeclareSymbol{IXXiI'Xi_typed*}{-1}{\draw (-1,1) node[xix-green-red1] {} -- (0,-0.3);
\draw[kernels2,rednode] (1,1) node[xigreen1] {} -- (0,-0.3) {};}

\DeclareSymbol{XiX}{-2.8}{\node[xibx] {};}
\DeclareSymbol{XXi}{-2.8}{\node[xibx] {};}
\DeclareSymbol{tauX}{-2.8}{ \draw[kernels2] (0,0) node[xibx] {};}
\DeclareSymbol{Xi}{-2.8}{\node[xib] {};}
\DeclareSymbol{Xigreen}{-2.8}{\node[xigreen] {};}
\DeclareSymbol{Xired}{-2.8}{\node[xired] {};}
\DeclareSymbol{not}{-2.8}{\node[not,minimum size=1.6mm] {};}
\DeclareSymbol{XiR}{-2.8}{\node[xie,fill=rednode] {};}
\DeclareSymbol{XiY}{-2.8}{\node[xie] {};}
\DeclareSymbol{XiZ}{-2.8}{\node[xid] {};}
\DeclareSymbol{IXiX}{0}{\draw (0,-0.25) node[not] {} -- (0,1.5) node[xix] {};}

\makeatletter 
\newcommand*{\bigcdot}{}
\DeclareRobustCommand*{\bigcdot}{%
  \mathbin{\mathpalette\bigcdot@{}}%
}
\newcommand*{\bigcdot@scalefactor}{.5}
\newcommand*{\bigcdot@widthfactor}{1.15}
\newcommand*{\bigcdot@}[2]{%
  \sbox0{$#1\vcenter{}$}
  \sbox2{$#1\cdot\m@th$}%
  \hbox to \bigcdot@widthfactor\wd2{%
    \hfil
    \raise\ht0\hbox{%
      \scalebox{\bigcdot@scalefactor}{%
        \lower\ht0\hbox{$#1\bullet\m@th$}%
      }%
    }%
    \hfil
  }%
}
\makeatother

\newcommand{\cut}{\mathfrak{C}}

\newcommand{\mrd}{\mathop{}\!\mathrm{d}}

\newcommand{\mcC}{\mathcal{C}}
\newcommand{\mcB}{\mathcal{B}}

\newcommand{\mcI}{\mathcal{I}}
\newcommand{\mcF}{\mathcal{F}}

\newcommand{\mcD}{\mathcal{D}}

\newcommand{\mcP}{\mathcal{P}}

\newcommand{\mcX}{\mathcal{X}}

\newcommand{\mcZ}{\mathcal{Z}}

\newcommand{\mbi}{\mathbf{i}}

\newcommand{\T}{\mathbf{T}}

\def\${|\!|\!|}

\def\gr#1{#1\textnormal{-gr}}

\def\axis{\textnormal{ax}}

\def\hol{\textnormal{hol}}

\def\scal#1{{\langle#1\rangle}}

\def\sol{{\mathop{\mathrm{sol}}}}

\newcommand{\mfu}{\mathfrak{u}}

\newcommand{\mfR}{\mathfrak{R}}

\newcommand{\mfp}{\mathfrak{p}}

\newcommand{\mfG}{\mathfrak{G}}
\newcommand{\mfO}{\mathfrak{O}}

\newcommand{\mfg}{\mathfrak{g}}

\def\cC{\mathscr{C}}

\def\SYM{\textnormal{\scriptsize \textsc{sym}}}

\def\BPHZ{\textnormal{\scriptsize \textsc{bphz}}}

\newcommand{\ad}{\mathrm{ad}}
\newcommand{\Ad}{\mathrm{Ad}}

\newcommand{\Trace}{\mathrm{Tr}}

\def\combplus[#1,#2,#3,#4]{\binom{#1\ {\scriptstyle #4} }{#2\ #3}}

\def\singlescalegenvert[#1,#2]{\hat{H}^{#2}_{#1}}
\def\multiscalegenvert[#1,#2]{H^{#2}_{#1}}

\def\moll{\chi}

\def\nr[#1]{\tilde{N}[#1]} 
\def\inn[#1]{\mathring{N}[#1]}
\def\nrinn[#1]{\hat{N}_{#1}} 
\def\nrmod[#1,#2]{\tilde{N}_{#1}(#2)}
\def\nrinnmod[#1,#2]{\hat{N}_{#1}(#2)}

\def\ident[#1]{\underline{#1}}

\def\mylink#1#2{\mathrel{\vbox{\offinterlineskip\ialign{%
    \hfil##\hfil\cr
    $\scriptscriptstyle#1$\cr
    \noalign{\kern0.1ex}
    $#2$\cr
}}}}
\def\mysublink[#1]#2#3{\mathrel{\vbox{\offinterlineskip\ialign{%
    \hfil##\hfil\cr
    $\scriptscriptstyle#2$\cr
    \noalign{\kern0.1ex}
    $#3$\cr
    \noalign{\kern-0.2ex}
    \smash{\raisebox{-\height}{\hbox{$\scriptscriptstyle #1$}}}\cr
    \noalign{\kern0.2ex}
}}}}

\def\fon[#1]{\cC_{#1}}

\def\mincompproj[#1]{\mfp_{#1}}

\def\Proj_#1{\mathop{\mathrm{Proj}_{#1}}}

\def\negrenorm[#1]{\mfR_{#1}}
\def\topnegrenorm[#1]{\overline{\mfR}_{#1}}

\def\quotedge[#1]{E^{q}_{#1}}

\def\posrenorm[#1]{\mcC_{#1}}
\def\topposrenorm[#1]{\overline{\mcC_{#1}}}
\def\cutsmod[#1]{\mathbb{C}_{+,#1}}

\def\fullcutsmod[#1]{\cut_{#1}}

\colorlet{symbols}{blue!30!black!50}
\colorlet{testcolor}{green!60!black}
\colorlet{darkblue}{blue!60!black}
\colorlet{darkgreen}{green!60!black}
\definecolor{darkergreen}{rgb}{0.0, 0.5, 0.0}

\definecolor{purple}{rgb}{0.55,0.05,0.8}

\colorlet{redkernel}{red!80}

\def\symbol#1{\textcolor{symbols}{#1}}

\def\1{\mathbf{\symbol{1}}}

\makeatletter
\pgfdeclareshape{crosscircle}
{
  \inheritsavedanchors[from=circle] 
  \inheritanchorborder[from=circle]
  \inheritanchor[from=circle]{north}
  \inheritanchor[from=circle]{north west}
  \inheritanchor[from=circle]{north east}
  \inheritanchor[from=circle]{center}
  \inheritanchor[from=circle]{west}
  \inheritanchor[from=circle]{east}
  \inheritanchor[from=circle]{mid}
  \inheritanchor[from=circle]{mid west}
  \inheritanchor[from=circle]{mid east}
  \inheritanchor[from=circle]{base}
  \inheritanchor[from=circle]{base west}
  \inheritanchor[from=circle]{base east}
  \inheritanchor[from=circle]{south}
  \inheritanchor[from=circle]{south west}
  \inheritanchor[from=circle]{south east}
  \inheritbackgroundpath[from=circle]
  \foregroundpath{
    \centerpoint%
    \pgf@xc=\pgf@x%
    \pgf@yc=\pgf@y%
    \pgfutil@tempdima=\radius%
    \pgfmathsetlength{\pgf@xb}{\pgfkeysvalueof{/pgf/outer xsep}}%
    \pgfmathsetlength{\pgf@yb}{\pgfkeysvalueof{/pgf/outer ysep}}%
    \ifdim\pgf@xb<\pgf@yb%
      \advance\pgfutil@tempdima by-\pgf@yb%
    \else%
      \advance\pgfutil@tempdima by-\pgf@xb%
    \fi%
    \pgfpathmoveto{\pgfpointadd{\pgfqpoint{\pgf@xc}{\pgf@yc}}{\pgfqpoint{-0.707107\pgfutil@tempdima}{-0.707107\pgfutil@tempdima}}}
    \pgfpathlineto{\pgfpointadd{\pgfqpoint{\pgf@xc}{\pgf@yc}}{\pgfqpoint{0.707107\pgfutil@tempdima}{0.707107\pgfutil@tempdima}}}
    \pgfpathmoveto{\pgfpointadd{\pgfqpoint{\pgf@xc}{\pgf@yc}}{\pgfqpoint{-0.707107\pgfutil@tempdima}{0.707107\pgfutil@tempdima}}}
    \pgfpathlineto{\pgfpointadd{\pgfqpoint{\pgf@xc}{\pgf@yc}}{\pgfqpoint{0.707107\pgfutil@tempdima}{-0.707107\pgfutil@tempdima}}}
  }
}
\makeatother

\colorlet{greennode}{green!50!black}
\colorlet{rednode}{red!50!black}
\colorlet{lbluenode}{blue!25}
\colorlet{dbluenode}{blue}
\colorlet{orangenode}{orange}

\definecolor{connection}{rgb}{0.7,0.1,0.1}

\tikzset{
root/.style={circle,fill=black!50,inner sep=0pt, minimum size=3mm},
        var/.style={circle,fill=black!10,draw=black,inner sep=0pt, minimum size=3mm},
        kernel/.style={semithick,shorten >=2pt,shorten <=2pt},
        kernel1/.style={thick},
        kernels/.style={snake=zigzag,shorten >=2pt,shorten <=2pt,segment amplitude=1pt,segment length=4pt,line before snake=2pt,line after snake=5pt,},
        rho/.style={densely dashed,semithick,shorten >=2pt,shorten <=2pt},
           testfcn/.style={dotted,semithick,shorten >=2pt,shorten <=2pt},
           tau/.style={circle,inner sep=1pt,draw=black,fill=white,text=black,thin},
        renorm/.style={shape=circle,fill=white,inner sep=1pt},
        labl/.style={shape=rectangle,fill=white,inner sep=1pt},
        xi/.style={very thin,circle,fill=lbluenode,draw=symbols,inner sep=0pt,minimum size=1.2mm},
        xigreen/.style={very thin,circle,fill=greennode,draw=symbols,inner sep=0pt,minimum size=1.2mm},
        xigreen1/.style={very thin,rectangle,fill=greennode,draw=symbols,inner sep=0pt,minimum size=1.2mm},
        xired/.style={very thin,circle,fill=rednode,draw=symbols,inner sep=0pt,minimum size=1.2mm},
        xilblue/.style={very thin,circle,fill=lbluenode,draw=symbols,inner sep=0pt,minimum size=1.2mm},
        xiorange/.style={very thin,circle,fill=orangenode,draw=symbols,inner sep=0pt,minimum size=1.2mm},
        xix/.style={crosscircle,fill=lbluenode,draw=symbols,inner sep=0pt,minimum size=1.2mm},
xix-green-red/.style={circle, fill=greennode!70!white,draw=rednode,inner sep=0pt,minimum size=1.6mm,append after command={node [inner sep=0pt,minimum size=0.8mm,thick, draw = rednode, cross out]{}}},
xix-green-red1/.style={rectangle, fill=greennode!70!white,draw=rednode,inner sep=0pt,minimum size=1.5mm,append after command={node [inner sep=0pt,minimum size=1mm,thick, draw = rednode, cross out]{}}},
	xib/.style={very thin,circle,fill=lbluenode,draw=symbols,inner sep=0pt,minimum size=1.6mm},
	xie/.style={very thin,circle,fill=greennode,draw=symbols,inner sep=0pt,minimum size=1.6mm},
	xid/.style={very thin,circle,fill=lbluenode,draw=symbols,inner sep=0pt,minimum size=1.6mm},
	xibx/.style={crosscircle,fill=lbluenode,draw=symbols,inner sep=0pt,minimum size=1.6mm},
	kernels2/.style={ultra thick,draw=symbols,segment length=12pt},
	not/.style={thin,regular polygon, regular polygon sides=3,draw=connection,fill=connection,inner sep=0pt,minimum size=1.2mm},
	notlblue/.style={thin,regular polygon, regular polygon sides=3,draw=lbluenode,fill=lbluenode,inner sep=0pt,minimum size=1.2mm},
	notorange/.style={thin,regular polygon, regular polygon sides=3,draw=orangenode,fill=orangenode,inner sep=0pt,minimum size=1.2mm},
	notgreen/.style={thin,regular polygon, regular polygon sides=3,draw=greennode,fill=greennode,inner sep=0pt,minimum size=1.2mm},
	>=stealth,
  }

\colorlet{darkblue}{blue!90!black}
\colorlet{darkred}{red!90!black}
\colorlet{darkgreen}{green!70!black}

\def\${|\!|\!|}

\def\?{{\color{red}?}}

\def\restr{\mathbin{\upharpoonright}}

\def\Gauge{\textnormal{\tiny {Gauge}}}

\def\dash{\leavevmode\unskip\kern0.18em--\penalty\exhyphenpenalty\kern0.18em}
\def\slash{\leavevmode\unskip\kern0.15em/\penalty\exhyphenpenalty\kern0.15em}

\let\basepoint\logof
\def\logof{\mathord{{\basepoint}}} 

\newcommand{\init}{\mathcal{I}}
\newcommand{\state}{\mathcal{S}}

\usepackage{stmaryrd}

\def\heatgr#1{{|\!|\!| #1 |\!|\!|}}

\numberwithin{equation}{section}

\let\OLDthebibliography\thebibliography
\renewcommand\thebibliography[1]{
  \OLDthebibliography{#1}
  \setlength{\parskip}{0pt}
  \setlength{\itemsep}{0.6mm}
}

\title{Stochastic quantisation of Yang--Mills}

\author{
Ilya Chevyrev
\thanks{School of Mathematics, The University of Edinburgh, Edinburgh EH9 3FD, United Kingdom.
\href{mailto:ichevyrev@gmail.com}{\tt ichevyrev@gmail.com}}
}

\date{}

\begin{document}

\maketitle

\begin{abstract}
We review two works~\cite{CCHS20,CCHS22} which study the stochastic quantisation equations of Yang--Mills on two and three dimensional Euclidean space with finite volume. The main result of these works is that one can renormalise the 2D and 3D stochastic Yang--Mills heat flow so that the dynamic becomes gauge covariant in law. Furthermore, there is a state space of distributional $1$-forms $\mathcal{S}$ to which gauge equivalence $\sim$ extends and such that the renormalised stochastic Yang--Mills heat flow projects to a Markov process on the quotient space of gauge orbits $\mathcal{S}/{\sim}$. In this review, we give unified statements of the main results of these works, highlight differences in the methods, and point out a number of open problems.
\end{abstract}

\renewcommand{\thefootnote}{\fnsymbol{footnote}} 
\renewcommand{\thefootnote}{\arabic{footnote}} 

\tableofcontents

\section{Introduction}

\subsection{Yang--Mills theory}

Yang--Mills (YM) theory plays an important role in the description of force-carrying particles in the Standard model.
An important unsolved problem in mathematics
is to show that YM theory on Minkowski space-time can be rigorously quantised.
We refer to~\cite{JW06} for a description of this problem, together with the surveys~\cite{Chatterjee18,LevySengupta17} for
related literature and problems.

A basic ingredient in YM theory
is a compact Lie group $G$, called the \textit{structure group}, with Lie algebra $\mfg$.
Passing to Euclidean space and working in an arbitrary dimension $d$, one can reformulate the problem as trying to make sense of the
YM probability measure on the space of $\mfg$-valued $1$-forms $A=(A_1,\ldots, A_4)\colon \R^d\to \mfg^d$
\[
\mu(\mcD A) = \mcZ^{-1} \exp(-S(A)) \mcD A\;,
\]
where $\mcD A$ on the right-hand side is a formal Lebesgue measure, $\mcZ$ is a normalisation constant, and $S(A)$ is the \textit{YM action}
\begin{equation*}
S(A) = \int_{\R^d} |F(A)|^2\;.
\end{equation*}
Above, $F(A)$ is the curvature $2$-form of $A$
given by $F_{ij}(A) = \partial_i A_j - \partial_j A_i + [A_i,A_j]$,
and we equip $\mfg$ with an $\Ad$-invariant inner product $\scal{\cdot,\cdot}$ and norm $|\cdot|$.

\begin{remark}
Without loss of generality, we can take $G\subset U(N)$ and $\mfg\subset \mfu(N)$ for some $N\geq 1$.
In this case, the adjoint action is $\Ad_gA = gAg^{-1}$
and a possible choice for the inner product is $\scal{X,Y} = -\Trace(XY)$,
and one can rewrite $|F(A)|^2 = -\frac12\sum_{i,j=1}^d \Trace(F_{ij}(A)^2)$.
\end{remark}

The case $d=4$ corresponds to physical space-time, but the task of constructing the probability measure $\mu$ makes sense for arbitrary dimension,
and even for a Riemannian manifold $M$ in place of $\R^d$.
The $\mfg$-valued $1$-forms in this case become connections on a principal $G$-bundle $P\to M$,
and one aims to define the probability measure $\mu$ on the space of connections on $P$.

The cases $d=2,3$ are considered substantially simpler than $d=4$ as they correspond to super-renormalisable
theories in quantum field theory (vs. renormalisable for $d=4$ and non-renormalisable for $d\geq 5$).
In the remainder of the article, we will focus on these dimensions and further restrict to \textit{finite volume}, replacing $\R^d$ by the torus $\T^d=\R^d/\Z^d$.
The underlying principal bundle $P$ is always assumed trivial
and we keep in mind that all geometric objects (connections, curvature forms, etc.) can be written in coordinates 
($\mfg$-valued $1$-forms, $2$-forms, etc.) once we fix a global section of $P$ which identifies it with $\T^d\times G$.
The space of connections is an affine space, with the difference of two connections being a $1$-form.

An important postulate in the physics of YM theory (and gauge theories in general)
is that all physical quantities should be invariant under the action of the \textit{gauge group}, i.e.\ the automorphism group of $P$.
Since $P$ is trivial, the gauge group consists of functions $g\colon\T^d\to G$.
These gauge transformation change the global section we use to identify $P\simeq \T^d\times G$ and
act on $1$-forms by
\begin{equ}
A\mapsto A^g \eqdef \Ad_g A - (\mrd g) g^{-1}\;.
\end{equ}
One can view $A^g$ either as a new connection, gauge equivalent to $A$,
or simply as the same connection $A$ written in a new coordinate system.


We write $A\sim B$ if there exists $g$ such that $A^g=B$ and write $[A]=\{B\,:\,B\sim A\}$ for the \textit{gauge orbit} of $A$.
In light of the above, the natural space on which to define the probability measure $\mu$ is not the space of $1$-forms,
but rather the quotient space $\mfO$ of all gauge orbits.
The space $\mfO$ is a non-linear space if $G$ is non-Abelian, which makes non-trivial
even the construction of the state space on which the YM measure $\mu$ should be defined.

A number of works have made contributions to a precise definition of this measure.
The most successful case is dimension $d=2$,
which includes $\R^2$ and compact orientable surfaces.
The key feature which makes 2D YM special is its exact solvability,
which allows one to write down an explicit formula for the joint distribution of Wilson loop observables;
this property was observed in the physics literature by Migdal~\cite{Migdal75} and later developed in mathematics; see e.g.~\cite{GKS89, Driver89, Fine91, Sengupta97, Levy03, Levy06, Chevyrev19YM}.

In the Abelian case $G=U(1)$ on $\R^d$ one can make sense of the measure $\mu$ for $d=3$~\cite{Gross83} and $d=4$~\cite{Driver87}.
For any structure group $G$, a form of ultraviolet stability on $\T^d$ was demonstrated for $d=4$ using a continuum regularisation in~\cite{MRS93}
and for $d=3,4$ using a renormalisation group approach on the lattice in a series of works by Balaban~\cite{Balaban85IV,Balaban87,Balaban89} (see also~\cite{Federbush86}).
However, a construction of the 3D YM measure and a description of its gauge-invariant observables, even on $\T^3$,
remains open.

\subsection{Stochastic quantisation}

Another approach to the construction of the YM measure was recently initiated in~\cite{CCHS20,CCHS22}, which is based on stochastic quantisation
(see also~\cite{Shen21} which treats scalar QED on $\T^2$).
The basic idea behind this approach is to view the measure $\mu$ as the invariant measure of a Langevin dynamic.
By studying this dynamic, one can try to determine properties and even give constructions of $\mu$.
The method was put forward in the context of gauge theories by Parisi--Wu~\cite{ParisiWu},
and has recently been applied to the construction of \textit{scalar} theories, see~\cite{AK20,MoinatWeber20,GH21,HS22}.


The Langevin dynamic associated to the YM measure is
\begin{equ}\label{eq:SYM_unrenorm}
\partial_t A = -\mrd_A^* F(A) + \xi\;,
\end{equ}
where $\mrd^*_A$ is the adjoint of the covariant derivative $\mrd_A$
and $\xi$ is a white noise built over $L^2$-space of $\mfg$-valued $1$-forms.

We point out right away that a difficulty in solving~\eqref{eq:SYM_unrenorm},
even in the absence of noise, is that the equation $\partial_t A = -\mrd_A^* F(A)$ is not parabolic.
This is a well-known feature of YM theory and is connected with the infinite-dimensional nature of the gauge group:
the YM equations $\mrd^*_AF(A) = 0$ are not elliptic and admit infinitely many solutions since, if $A$ is a solution, then so is every element of the gauge orbit $[A]$.

Fortunately, it is possible to bypass this issue by adding a gauge-breaking term to the right-hand side of~\eqref{eq:SYM_unrenorm}
which is tangent to the orbit $[A]$ at $A$ and which renders the equation parabolic.
This `trick' was used by Zwanziger~\cite{zwanziger81} and Donaldson~\cite{Donaldson}
in YM theory (see also~\cite{BHST87II,Sadun}), and by DeTurck~\cite{deturck83} in the context of Ricci flow.
The tangent space of $[A]$ at $A$ is precisely the set of all $\mrd_A\omega=\sum_{i=1}^d( \partial_i + [A_i,\cdot])\omega \mrd x_i$
where $\omega$ is a $\mfg$-valued function.
A natural choice for the additional term is $-\mrd_A\mrd^* A$, where $\mrd^* A = -\sum_{i=1}^d \partial_iA_i$.
Therefore, making sense of~\eqref{eq:SYM_unrenorm} is equivalent to making sense of the parabolic equation
\begin{equation}\label{eq:SYM_DeTurck_no_renorm}
\partial_t A = -\mrd^*_A F(A) - \mrd_A \mrd^* A + \xi\;,
\end{equation}
which is the equation we consider henceforth.
In coordinates,~\eqref{eq:SYM_DeTurck_no_renorm} reads
\begin{equ}\label{eq:SYM_coord}
\partial_t A_i = \Delta A_i + [A_j,2\partial_jA_i - \partial_iA_j+[A_j,A_i]] + \xi_i\;,\quad i=1,\ldots, d\;,
\end{equ}
where $\xi_1,\ldots,\xi_d$ are i.i.d. $\mfg$-valued white noises and $[\cdot,\cdot]$ is the Lie bracket of $\mfg$. Here and below, there is an implicit summation over $j=1,\ldots, d$.

\subsubsection{Gauge covariance}
\label{subsec:gauge_covar}

The reason why~\eqref{eq:SYM_coord} is natural
is that it is (formally) gauge covariant in law.
To see this,
it is convenient to work in coordinate-free notation and, for the moment, make an distinction between connections and $1$-forms.
Recall that the space of connections is affine modelled on the space $\Omega^1$ of $\mfg$-valued $1$-forms (we are using the global section of $P$ to identify $\ad(P)$-valued forms with $\mfg$-valued forms).
For a connection $A$, recall further that $\mrd_A\colon\Omega^k\to \Omega^{k+1}$ is a linear map from $\mfg$-valued $k$-forms to $(k+1)$-forms
with adjoint $\mrd^*_A\colon \Omega^{k+1}\to\Omega^k$.
Furthermore, gauge transformations $g$ act on connections via $A\mapsto A^g$ and on forms via $\omega\mapsto \Ad_g \omega$.

Consider now a time-dependent connection $A$ and gauge transformation $g$.
Then $B \eqdef A^g$ satisfies
\[
\partial_t B
=
 \Ad_g(\partial_t A)
- \mrd_B[(\partial_t g) g^{-1}]\;.
\]
Let $Z$ denote the canonical flat connection associated with our choice of global section, i.e.\ the connection associated with the $1$-form $0$.
Consider a time-dependent $1$-form $\xi$ and
suppose that $A$ solves~\eqref{eq:SYM_coord},
which we now write as
\begin{equ}
\partial_t A = -\mrd^*_A F(A) - \mrd_A \mrd^*_A (A-Z) + \xi \;.
\end{equ}
Observe  the covariance properties
\begin{equ}
\Ad_g(A-Z)=A^g-Z^g\;, \quad
\Ad_g (\mrd_A \omega) = \mrd_{A^g}(\Ad_g \omega)\;,
\quad \Ad_g [F(A)] = F(A^g)
\end{equ}
(the whole reason for the above notation is that these properties take on a simple form).
It follows that
$B$ solves
\[
\partial_t B = -\mrd^*_B F(B) - \mrd_B \mrd^*_{B} (B-Z^g) + \Ad_g\xi
- \mrd_B[(\partial_t g) g^{-1}]\;.
\]
To bring $B$ into a form similar to $A$, it thus natural to take $g$ which solves
\begin{equ}
(\partial_t g) g^{-1} = \mrd^*_{B} (Z^g -Z)\;.
\end{equ}
In this case, the equation for $B$ becomes
\begin{equ}
\partial_t B = -\mrd^*_B F(B) - \mrd_B \mrd^*_{B} (B-Z) + \Ad_g\xi\;,
\end{equ}
which is almost the same as the equation for $A$, except that $\xi$ is replaced by $\Ad_g\xi$.

If we now assume $\xi$ is a white noise, and that the equations above make sense with global in time solutions,
then $\Ad_g \xi$ is equal in law to $\xi$ by It{\^o} isometry.
This formal argument suggests there is a coupling between two solutions $A,B$ to~\eqref{eq:SYM_coord} started from gauge equivalent initial conditions such that $A(t)\sim B(t)$ for all $t\geq 0$.
In particular, the law of the projected process $[A]$ on gauge orbits is equal to that of $[B]$.
The projected process on gauge orbits is therefore well-defined and Markov,
and its invariant probability measure is a natural candidate for the YM measure.

\subsection{Main results}
\label{subsec:main_results}

The basic objective of the works~\cite{CCHS20} and~\cite{CCHS22}
is to make rigorous the above formal argument in the case of $\T^d$ 
for $d=2$ and $d=3$ respectively.
In particular, they aim to define a natural Markov process on gauge orbits associated to the YM Langevin dynamic~\eqref{eq:SYM_coord}.
In this subsection, we give unified statements of the main results in~\cite{CCHS20,CCHS22}.
We will go into more detail and discuss the differences in proofs in Sections~\ref{sec:2D} and~\ref{sec:3D}.
We also discuss in detail these results in the simple case that $G=U(1)$ in Section~\ref{subsec:Abelian}.

\begin{remark}
In~\cite{CCHS22}, the more general Yang--Mills--Higgs theory is considered.
To simplify our discussion, we restrict here the results of~\cite{CCHS22} to pure YM theory.
\end{remark}

One of the main difficulties in solving~\eqref{eq:SYM_coord}
is that,
unless $G$ is Abelian, this equation is non-linear and \textit{singular}.
To simplify notation and highlight the nature of the non-linearities, we will henceforth write~\eqref{eq:SYM_DeTurck_no_renorm} and~\eqref{eq:SYM_coord} as
\begin{equ}\label{eq:SYM_simple}
\partial_t A = \Delta A + A\partial A + A^3 + \xi\;.
\end{equ}
It is well-known that white noise on $\R\times \T^d$ can be realised as a random distribution in
$C^{-1-\frac{d}{2}-\kappa}$, the H{\"o}lder--Besov space with parabolic scaling,
for arbitrary $\kappa>0$.
Furthermore, this regularity is optimal, at least in the scale of Besov spaces.
By Schauder estimates, we therefore expect that $A$ is in $C^{1-\frac d2 -\kappa}$ (and no better), but this renders the non-linear terms $A\partial A$ and $A^3$ analytically ill-defined once $d\geq 2$ because $1-\frac d2 -\kappa<0$.
(We recall here that the bilinear map $(f,g)\mapsto fg$ defined for smooth functions extends to $C^\alpha\times C^\beta$ if and only if $\alpha+\beta>0$.)

A number of solution theories, including regularity structures~\cite{Hairer14} and paracontrolled calculus~\cite{GIP15} (see also~\cite{Kupiainen16,OW19,Duch21}),
have been developed in the last decade which allow one to make sense of such equations.
The key condition which must be satisfied is \textit{subcriticality}, which happens if and only if $d<4$ and which parallels the notion of super-renormalisability in QFT.
Subcriticality implies that the solution $A$ to~\eqref{eq:SYM_simple} should be a perturbation of the stochastic heat equation (SHE)
\begin{equ}
\partial_t \Psi = \Delta \Psi + \xi\;,
\end{equ}
which has the Gaussian free field (GFF) as an invariant measure.
The solution theory which~\cite{CCHS20,CCHS22} use and build on is the theory of regularity structures.

Recall that $A$ represents a geometric object (a principal $G$-connection).
On the other hand, the solution $A$ to~\eqref{eq:SYM_simple} at positive times is expected to 
be a distribution of the same regularity as the GFF.
Therefore, a first natural question is whether there exists a state space $\state$ large enough to support the GFF while small enough so that gauge equivalence extends to $\state$.
One of the main results of~\cite{CCHS20,CCHS22}
gives an answer to this question, which can be informally stated as follows.

\begin{theorem}\label{thm:state_space}
There exists a metric space $(\state,\Sigma)$ of $\mfg$-valued distributional $1$-forms on $\T^d$, $d=2,3$,
which contains all smooth $1$-forms and to
which gauge equivalence $\sim$ extends in a canonical way.
Furthermore, $\state$ contains distributions of the same regularity as the GFF on $\T^d$.
\end{theorem}

The constructions of $\state$ in~\cite{CCHS20} and in~\cite{CCHS22} are rather different.
In~\cite{CCHS20}, $\state$ (therein denoted by $\Omega^1_\alpha$) is a Banach space defined through line integrals, and gauge equivalence is determined by an action of a gauge group.
In~\cite{CCHS22}, $\state$ is a non-linear metric space of distributions defined in terms of the effect of the heat flow; gauge equivalence is extended using a gauge-covariant regularising operator (the deterministic YM flow).
We describe these constructions further in Sections~\ref{subsec:2D_state} and~\ref{subsec:3D_state}.

\begin{remark}
A naive construction of $\state$ would be to take $C^\eta$ for $\eta= 1-\frac d2-\kappa$
and quotient by the action of the gauge group $C^{\eta+1}(\T^d,G) \cap C^{-\eta+\kappa}(\T^d,G)$,
so that $\Ad_g A\in C^\eta$ and
$(\mrd g) g^{-1} \in C^\eta$ are well-defined.
Such a construction ends up lacking most of the nice properties we discuss in Sections~\ref{subsec:2D_state} and~\ref{subsec:3D_state}.
\end{remark}

\begin{remark}\label{rem:quotient_top}
While Theorem~\ref{thm:state_space}
makes it seems like the 2D and 3D cases are on an equal footing, we actually know much more about $\state$ in 2D than in 3D,
e.g.\ the space of orbits $\state/{\sim}$ in 2D comes with a natural complete metric, and is thus Polish,
while we only know that $\state/{\sim}$ is completely Hausdorff (and separable) in 3D.
\end{remark}

We now turn to the question of solving~\eqref{eq:SYM_simple}.
A natural approach is to replace $\xi$ by a smooth approximation $\xi^\eps$ and let the mollification parameter $\eps\downarrow0$.
The hope then is that the corresponding solutions $A$ converge.
Unfortunately, this is not in general the case and the SPDE requires \textit{renormalisation}
for convergence to take place.
The following result ensures that renormalised solutions to~\eqref{eq:SYM_simple} exist, at least up to a potential finite time blow-up.

A \textit{mollifier} is a smooth compactly supported function $\moll \colon \R\times \R^d \to \R$ such that 
$\int\moll=1$
and which is spatially symmetric and invariant under flipping coordinates $x_i\mapsto-x_i$ for $i=1,\ldots, d$.
Denoting $\moll^\eps(t,x) = \eps^{-2-d}\moll(\eps^{-2}t,\eps^{-1}x)$
and $*$ for space-time convolution,
we define $\xi^\eps=\moll^\eps*\xi$.
Furthermore, define
\begin{equ}
L_G(\mfg,\mfg) = \{X\in L(\mfg,\mfg)\,:\, X\Ad_g = \Ad_gX \text{ for all } g\in G\}\;.
\end{equ}

\begin{theorem}\label{thm:local}
For every mollifier $\moll$, there exists a family of operators $\{C^\eps_{\BPHZ}\}_{\eps\in(0,1)}\subset L_G(\mfg,\mfg)$
such that, for any $\mathring C\in L(\mfg,\mfg)$ and initial condition $A(0)\in\state$,
the solution to the PDE
\begin{equ}\label{SYM}
\partial_t A = \Delta A + A\partial A + A^3 + \xi^\eps + (C^\eps_\BPHZ + \mathring C) A\;,
\tag{SYM}
\end{equ}
converges in probability in $C(\R_+,\state\sqcup\{\skull\})$ as $\eps\downarrow0$.
The limit as $\eps\downarrow0$ furthermore does not depend on $\moll$.
\end{theorem}

In coordinates,~\eqref{SYM} reads
\begin{equ}
\partial_t A_i = \Delta A_i + [A_j,2\partial_jA_i - \partial_iA_j+[A_j,A_i]] + \xi_i^\eps + (C^\eps_\BPHZ + \mathring C)A_i\;,\quad i=1,\ldots, d\;.
\end{equ}
Remark that the operator $C^\eps_\BPHZ + \mathring C \in L(\mfg,\mfg)$ in~\eqref{SYM} is `block diagonal' meaning that only $(C^\eps_\BPHZ + \mathring C)A_i$ appears in the equation for $A_i$.

\begin{definition}\label{def:SYM}
We call the $\eps\downarrow0$ limit of $A$ as in Theorem~\ref{thm:local} the \textit{solution to~\eqref{SYM} driven by $\xi$ with bare mass $\mathring C$.}
\end{definition}

The bare mass $\mathring C$ is used to parametrise the space of all `reasonable' solutions and is a free parameter at this stage.
We will see below (Theorems~\ref{thm:generative} and~\ref{thm:B_bar_A})
that there does exist a unique choice for $\mathring C$ which selects a distinguished element of this solution space.

\begin{remark}\label{rem:BPHZ_not_unique}
The operators $C^\eps_\BPHZ$ are called the \textit{BPHZ constants}
and are given by (in principle explicit) integrals involving $\moll$ and an arbitrary large scale truncation $K$ of the heat kernel.
While the solution to~\eqref{SYM} in Definition~\ref{def:SYM} is independent of $\moll$, it does in general depend on the choice of $K$ used to define $C^\eps_\BPHZ$.
\end{remark}

\begin{remark}
The point $\skull$ is a cemetery state and is added to $\state$ to handle the possibility of finite time blow-up.
Some care is needed to properly define $C(\R_+,\state\sqcup\{\skull\})$
and the metric that one equips it with.
This is done in~\cite[Sec.~1.5.1]{CCHS20},
where, for a general metric space $E$, a metric space $E^\sol$ of continuous paths with values in $E\sqcup \{\skull\}$
is defined in which two paths are close if they track each other until the point when they become large.
Our notation
$C(\R_+,\state\sqcup\{\skull\})$ here really means $\state^\sol$.
\end{remark}

\begin{remark}\label{rem:2D_const_conv}
It turns out that in 2D, due to a cancellation in renormalisation constants,
$C^\eps_\BPHZ$ converges to a finite value as $\eps\downarrow0$; see Theorem~\ref{thm:2D_local}.
Therefore, Theorem~\ref{thm:local} in 2D remains true if $C^\eps_\BPHZ+\mathring C$ replaced by any fixed $C\in L(\mfg,\mfg)$, which is the formulation of~\cite[Thm.~2.4]{CCHS20}.
No such cancellation occurs in 3D and
$C^\eps_\BPHZ$ diverges at order $\eps^{-1}$.
\end{remark}

We now discuss the way in which solutions to~\eqref{SYM} are gauge covariant in the sense described in Section~\ref{subsec:gauge_covar}.
Remark that, by inserting the counterterm $(C^\eps_\BPHZ + \mathring C) A$
we are seemingly breaking the desired gauge covariance property discussed in Section~\ref{subsec:gauge_covar} 
(in the notation of that section, $(C^\eps_\BPHZ + \mathring C) A$ should be written $(C^\eps_\BPHZ + \mathring C) (A-Z)$).
However, the formal argument in Section~\ref{subsec:gauge_covar} also breaks if we replace $\xi$ by $\xi^\eps$ because It\^o isometry is not true for the latter.

A surprising fact is that, if one chooses $\mathring C$ carefully, then the broken gauge covariance (due to the counterterm $(C^\eps_\BPHZ + \mathring C) A$) compensates in the $\eps\downarrow0$ limit for the broken It\^o isometry (due to the mollified noise $\xi^\eps$),
and one obtains a solution to~\eqref{SYM} which is gauge covariant in law.

It is not entirely trivial to make this statement precise, essentially because we do not know if~\eqref{SYM} (with any bare mass)
has global in time solutions.
In particular, we do not know how to rule out that solutions to~\eqref{SYM} with different gauge equivalent initial conditions $a\sim b$ blow up at different times,
and this makes it unclear in what sense we can expect the projected process $[A]$ on gauge orbits to be Markov.

To address this issue, it is natural to look for a type of process which solves~\eqref{SYM} on disjoint intervals $[\varsigma_{j-1},\varsigma_{j})$ and
at time $\varsigma_j$ jumps to a new representative of the gauge orbit $[\lim_{t\uparrow\varsigma_j}A(t)]$.
This should happen 
in such a way that $A$ does not blow  up unless the entire orbit $[A]$ `blows up'.
This class of processes
is defined through \textit{generative probability measures} in~\cite{CCHS20,CCHS22}.

We say that a probability measure $\mu$ on the space of c{\`a}dl{\`a}g functions $D(\R_+,\state\sqcup\{\skull\})$
is \textit{generative} with bare mass $\mathring C$ and initial condition $a\in \state$
if there exists a white noise $\xi$
and a random variable $A$ with law $\mu$
such that
\begin{enumerate}[label=(\roman*)]
\item\label{pt:init_cond} $A(0)=a$ almost surely,
\item\label{pt:SYM_jumps} there exists a sequence of stopping times $\varsigma_0=0\leq \varsigma_1\leq \varsigma_2<\ldots$ such that $A$ solves~\eqref{SYM} 
driven by $\xi$ with bare mass $\mathring C$
on each interval $[\varsigma_j,\varsigma_{j+1})$,

\item\label{pt:jumps_gauge} for every $j\geq 0$, $A(\varsigma_{j+1}) \sim \lim_{t\uparrow \varsigma_{j+1}} A(t)$, and

\item\label{pt:honest_blow_up} $\lim_{j\to\infty} \varsigma_j = T^* \eqdef \inf\{t\geq 0\,:\, A(t)=\skull\}$ and,
on the event $T^*<\infty$,\footnote{In 3D, one needs to work with a slightly stronger metric $\bar\Sigma$ in~\eqref{eq:2D_honest_blow_up}, see Section~\ref{subsec:3D_gauge_covar}.}
\begin{equ}\label{eq:2D_honest_blow_up}
\lim_{t\uparrow T^*} \inf_{B \sim A(t)} \Sigma(B,0) = \infty\;.
\end{equ}
\end{enumerate}


The point of this definition is to give a sufficiently general and natural way in which~\eqref{SYM} can be restarted along gauge orbits.
The following result from~\cite{CCHS20,CCHS22}
ensures the existence of a canonical Markov process associated to~\eqref{SYM} on the quotient space of gauge orbits $\mfO\eqdef \state/{\sim}$
\textit{provided} the bare mass is chosen in a precise way.

\begin{theorem}\label{thm:generative}
\begin{enumerate}[label=(\alph*)]
\item \label{pt:gen_exists} For every $a\in\state$ and $\mathring C\in L(\mfg,\mfg)$, there exists a generative probability measure $\mu$ with bare mass $\mathring{C}$ and initial condition $a$.

\item \label{pt:covar} There exists $\check C \in L_G(\mfg,\mfg)$ with the following properties.
For all $a\sim b\in\state$, if $\mu,\nu$ are generative probability measures with initial conditions $a,b$ respectively and  bare mass $\check C$,
then the pushforward measures $\pi_*\mu$ and $\pi_*\nu$ are equal.
In particular, the probability measure $\P^{x} = \pi_*\mu$, where $\mu$ is generative with bare mass $\check C$ and initial condition $a\in x\in\mfO$, depends only on $x$.
Finally, $\{\P^x\}_{x\in\mfO}$ are the transition functions of a time homogenous, continuous Markov process on $\mfO \sqcup\{\skull\}$.
\end{enumerate}
\end{theorem}

\begin{remark}
Theorem~\ref{thm:generative}\ref{pt:covar} makes no claims about
the \textit{uniqueness} of $\check C$, but we conjecture that $\check C$ is indeed unique (which is not difficult to prove in the Abelian case, see Section~~\ref{subsec:Abelian}).
\end{remark}

Finally, we mention a result in~\cite{CCHS20,CCHS22} which is crucial for the proof of Theorem~\ref{thm:generative}\ref{pt:covar}
and which makes precise the coupling argument outlined in Section~\ref{subsec:gauge_covar}.
This result makes a stronger statement about the constant $\check C$ for which uniqueness does hold.
It can also be seen
as a version of the Slavnov--Taylor identities
for renormalisation schemes that preserve gauge symmetries.

Suppose that $A$ solves~\eqref{SYM} and, recalling that $\Delta A + A\partial A + A^3 + CA$ is shorthand for $-\mrd^*_A F(A) - \mrd_A\mrd^*_A(A-Z) + C(A-Z)$,
suppose that $(B,g)$ solves
\begin{equs}[eq:B_no_coord]
\partial_t B &= -\mrd^*_B F(B) - \mrd_B \mrd^*_{B} (B-Z) + \Ad_g\xi^\eps
 + (C^\eps_\BPHZ + \mathring C) (B-Z^g)\;,
\\
(\partial_t g) g^{-1} &= \mrd_B^*(Z^g-Z)\;,
\end{equs}
where the initial condition of $B$ is $B(0)=A(0)^{g(0)}$.
Then, the same computation as in Section~\ref{subsec:gauge_covar} (see also~\cite[Sec.~2.2]{CCHS20}) shows that
$A^g = B$.
In coordinates,~\eqref{eq:B_no_coord}
is written as
\begin{equs}[eq:SPDE_B]
\partial_t B &= \Delta B + B\partial B + B^3
+\Ad_g \xi^\eps
+ ( C^\eps_\BPHZ + \mathring C)(B + (\mrd g) g^{-1})\;,
\\
(\partial_t g)g^{-1} &= \partial_j((\partial_j g )g^{-1}) + [B_j,(\partial_j g) g^{-1}]\;.
\end{equs}
It is now natural to compare~\eqref{eq:SPDE_B} to
\begin{equs}[eq:SPDE_bar_A_desired]
\partial_t \bar A &= \Delta \bar A + \bar A\partial \bar A + \bar A^3
+\moll^\eps*(\Ad_{\bar g} \xi)
+ ( C^\eps_\BPHZ + \mathring C)\bar A\;,
\\
(\partial_t \bar g)\bar g^{-1} &= \partial_j((\partial_j \bar g) \bar g^{-1}) + [\bar A_j,(\partial_j \bar g) \bar g^{-1}]\;,
\\
\bar g(0) &= g(0)\;,\quad \bar A(0) = B(0)\;.
\end{equs}
Remark that $\moll^\eps*(\Ad_{\bar g} \xi)$ is equal in law to $\xi^\eps$, and thus $\bar A$ is equal in law to the solution of~\eqref{SYM} with initial condition $B(0)$,
provided we take $\moll$ non-anticipative in the following sense.

\begin{definition}
A mollifier $\moll$ is called \textit{non-anticipative} if it has support in $(0,\infty)\times\R^d$.
\end{definition}

What we would therefore like to show is that, for a special choice of $\mathring C$,
the solutions to~\eqref{eq:SPDE_B}
and~\eqref{eq:SPDE_bar_A_desired} converge as $\eps\downarrow0$ to the same limit.
The identity $A^g=B$, which survives in the limit,
would provide a coupling between~\eqref{SYM} with initial condition $A(0)$ and initial condition $A(0)^{g(0)}$
under which the two solutions are gauge equivalent, at least locally in time.
It turns out that the following more general result is true.

\begin{theorem}\label{thm:B_bar_A}
There exists a unique
$\check C \in L_G(\mfg,\mfg)$
such that for all non-anticipative mollifiers $\moll$, all $\mathring C\in L(\mfg,\mfg)$,
and all initial conditions $(B(0),g(0))  \in \state\times C^\rho(\T^d,G)$, $\rho\in(\frac12,1)$,
the solution $(B,g)$ to~\eqref{eq:SPDE_B} 
converges as $\eps\downarrow0$ in probability to the same limit as the solution to
\begin{equs}[eq:SPDE_bar_A]
\partial_t \bar A &= \Delta \bar A + \bar A\partial \bar A + \bar A^3
+\moll^\eps*(\Ad_{\bar g} \xi)
+ ( C^\eps_\BPHZ + \mathring C)\bar A + (\mathring C - \check C)(\mrd \bar g)\bar g^{-1}\;,
\\
(\partial_t \bar g)\bar g^{-1} &= \partial_j((\partial_j \bar g) \bar g^{-1}) + [\bar A_j,(\partial_j \bar g) \bar g^{-1}]\;,
\\
\bar g(0) &= g(0)\;,\quad \bar A(0) = B(0) \;.
\end{equs}
Furthermore, the solution to~\eqref{SYM} with bare mass $\check C$ is independent of
$\moll$ and of the choice of $C^\eps_\BPHZ$.
\end{theorem}

It follows from Theorem~\ref{thm:B_bar_A}
that the solutions to~\eqref{eq:SPDE_B}
and~\eqref{eq:SPDE_bar_A_desired} indeed converge to the same limit as $\eps\downarrow0$
provided we choose $\mathring C= \check C$.
The operator $\check C$ in Theorem~\ref{thm:B_bar_A} is exactly the operator appearing in Theorem~\ref{thm:generative}\ref{pt:covar};
in 2D, we can give an explicit expression for $\check C$, see~\eqref{eq:check_C_2D}.

\begin{remark}
The value of $\check C$ in Theorem~\ref{thm:B_bar_A} is determined uniquely \textit{after} we fix a choice for $C^\eps_\BPHZ$.
However, recall from Remark~\ref{rem:BPHZ_not_unique} that $C^\eps_\BPHZ$ is not unique or canonical --- it is determined by $\moll$ and an arbitrary truncation of the heat kernel (see, e.g.\ Theorem~\ref{thm:2D_local}).
The final part of Theorem~\ref{thm:B_bar_A}
states that the solution of~\eqref{SYM} with bare mass $\check C$
is independent of $\moll$ and this choice of truncation.
\end{remark}

%

\subsection{Abelian case}
\label{subsec:Abelian}

We end this section by discussing the above results in the Abelian case, i.e.\ $G=U(1)$ and $\mfg=\R$.
We consider here $d\geq 1$ arbitrary.
We will see in this case that
\begin{itemize}
\item the constant in Theorem~\ref{thm:generative}\ref{pt:covar} is $\check C=0$ and is \textit{unique},
\item if $\T^d$ is replaced by $\R^d$, then uniqueness of $\check C$ fails (Remark~\ref{rem:Rd_nonunique}),
\item \eqref{SYM} with bare mass $\check C=0$ has global in time solutions but no
invariant probability measure
(Remark~\ref{rem:no_invar_measure})
\end{itemize}
In the Abelian case, the non-linearities $A\partial A$ and $A^3$ as well as the constants $C^\eps_\BPHZ$ vanish.
Equation~\eqref{SYM} with bare mass $\mathring C$ therefore becomes linear in $A$
and converges as $\eps\downarrow0$ to the solution of the SHE with a mass term
\begin{equ}\label{eq:SHE_Abelian}
\partial_t A = \Delta A + \mathring C A+\xi\;.
\end{equ}
Treating $1$-forms on $\T^d$ as periodic functions (or distributions) on $\R^d$ modulo $\Z^d$,
it is easy to see that
$A\sim B$ if and only if $A=B+\mrd \omega$ for some $\omega\colon \R^d\to \R$ such that $e^{\mbi \omega}$ is periodic where $\mbi=\sqrt{-1}$.
The tangent space of every gauge orbit is therefore $\{\mrd \omega\,:\, \omega\colon \T^d\to\mbi\R\}$.

Since $\Ad_g$ is now the identity, it is
clear that a possible value for $\check C$ in Theorems~\ref{thm:generative} and~\ref{thm:B_bar_A} is $\check C =0$.
This is because, if $B(0)=A(0)+\mrd \omega(0)$ and $A,B$ solve~\eqref{eq:SHE_Abelian} with $\mathring C=0$, then
$B = A+\mrd \omega$ for all times where $\omega$ solves $\partial_t \omega = \Delta \omega$.

Furthermore, $\check C=0$ is the \textit{only} possible value for $\check C$.
Indeed, consider the two gauge equivalent initial conditions $A(0)\eqdef 0$
and $B(0)=(B_1(0),\ldots,B_d(0))\eqdef (2\pi,0,\ldots, 0)$,
and suppose that $A,B$ solve~\eqref{eq:SHE_Abelian} with  $\mathring  C \neq 0$.
Then
$B(t) = e^{t\mathring C}B(0) + A(t)$,
where we used that $B(0)$ is constant on $\T^d$.
Consider now the gauge-invariant observable $I[A]
\eqdef
e^{\mbi \int_{\T^d} A_1}$.
Then
$\int_{\T^d} B_1(t) = e^{t\mathring C}2\pi +\int_{\T^d} A_1(t)$,
and thus
\begin{equ}
\E I[A(t)] = \exp(\mbi e^{t\mathring C}2\pi) \E I[B(t)] \neq \E I[B(t)]
\end{equ}
for all $t>0$ sufficiently small since $\mathring C\neq 0$.
This show that $A$ and $B$ cannot be gauge equivalent in law, and thus $\check C=0$ is the only possible value in Theorem~\ref{thm:generative}\ref{pt:covar}.

\begin{remark}\label{rem:Rd_nonunique}
The above argument relies on the fact that $\T^d$ is
not simply connected: we exploited
that $(\mrd g) g^{-1}$, which appears in~\eqref{eq:SPDE_bar_A}, is not tangent to gauge orbits in general (it is tangent if and only if $g=e^{\mbi\omega}$ for some $\omega\colon \T^d\to \R$).
If $\T^d$ is replaced by $\R^d$,
then $(\mrd g) g^{-1}$ \textit{is} tangent to gauge orbits
since we can always write $g=e^{\mbi\omega}$
for some $\omega\colon\R^d\to\R$.
Therefore, on $\R^d$ in the Abelian case, there is no uniqueness of $\check C$.
Explicitly, working on $\R^d$, suppose $A$ solves~\eqref{eq:SHE_Abelian}
and consider $B$ as in~\eqref{eq:SPDE_B}
where $B(0)=A(0)^{g(0)}$
but now $g$ satisfies
\begin{equ}
\mrd[(\partial_t g) g^{-1}] = -\mrd(\mrd^* [(\mrd g) g^{-1}]) + \mathring C(\mrd g) g^{-1}\;,
\end{equ}
i.e.\ $g=e^{\mbi\omega}$ where $\omega$ solves $
\partial_t \omega =  \Delta\omega+ \mathring C \omega$.
Then $B(t)=A(t)^{g(t)}= A(t) - \mrd \omega(t)$ for all $t\geq0$ and $B$ solves the same equation~\eqref{eq:SHE_Abelian} as $A$ for any $\mathring C\in\R$.
\end{remark}

\begin{remark}\label{rem:no_invar_measure}
Clearly~\eqref{eq:SHE_Abelian} with $\mathring C=0$ does not have an invariant \textit{probability} measure because $\int_{\T^d} A_i$ evolves like a Brownian motion.
In fact, any gauge equivalent generalisation of~\eqref{SYM} of the form $\partial_t A = \Delta A + \xi + \mrd\omega$, where $\omega$ is adapted,
will lack an invariant probability measure because the spatial mean is unaffected by $\mrd\omega$.
But we do obtain an invariant probability measure for the projected process $[A]$
because any $1$-form is gauge equivalent to another $1$-form $B$ such that $\int_{\T^d} B_1, \dots, \int_{\T^d}B_d \in [-\pi,\pi)$.
This remark shows that the Markov process from Theorem~\ref{thm:generative} can have an invariant probability measure while~\eqref{SYM} with bare mass $\check C$ does not.
\end{remark}

\section{Two dimensions}\label{sec:2D}

We describe in this section the main results in~\cite{CCHS20}, i.e.\ the results of Section~\ref{subsec:main_results}
in 2D.

\subsection{State space}
\label{subsec:2D_state}

The definition of the state space $\state$ (denoted by $\Omega^1_\alpha$ in~\cite{CCHS20}) is motivated by the desire to define holonomies, and thus Wilson loops, for every element $A\in\state$.
The construction is a refinement of that introduced in~\cite{Chevyrev19YM}.
Let $\mcX = \T^2 \times B_{1/4}$, where $B_{1/4} = \{v\in\R^2\,:\,|v|\leq \frac14\}$.
We think of $\mcX$ as the collection of straight line segments $\ell=(x,v)$ in $\T^2$ of length at most $|\ell|\eqdef |v|<\frac14$. (The starting point of $\ell$ is $x$.)

For $\alpha\in [0,1]$ and a smooth $1$-form $A\in C^\infty(\T^2,\mfg^2)$,
define the norm
\[
|A|_{\gr\alpha} = \sup_{\ell} \frac{|A(\ell)|}{|\ell|^\alpha}\;,
\]
where the supremum is taken over all $\ell\in\mcX$ with $|\ell|>0$ and where we define the \textit{line integral}
\begin{equ}\label{eq:line_int_def}
A(\ell) \eqdef \int_0^1 A(x+vt)v \mrd t = \int_0^1 \sum_{i=1}^2 A_i(x+vt)v_i \mrd t\;.
\end{equ}
We furthermore define the strengthened norm
\[
|A|_\alpha  = |A|_{\gr\alpha} + \sup_{P} \frac{|A(\partial P)|}{|P|^{\alpha/2}}\;,
\]
where the supremum is over all oriented triangles $P=(\ell_1,\ell_2,\ell_3)$ with $\ell_i\in \mcX$ and area $|P|>0$, and $A(\partial P) \eqdef \sum_{i=1}^3 A(\ell_i)$.
We can now define the state space studied in~\cite{CCHS20}.

\begin{definition}
The Banach space $(\state,|\cdot|_\alpha)$ is defined as the completion of smooth $\mfg$-valued $1$-forms under $|\cdot|_\alpha$ for some $\alpha\in(\frac23,1)$.
\end{definition}

\noindent
The metric $\Sigma$ in Theorem~\ref{thm:state_space}
is then the usual metric $\Sigma(A,B)=|A-B|_\alpha$.

\begin{remark}
To motivate these norms, consider $A=(A_1,A_2)$ a pair of i.i.d. GFFs.
A simple calculation shows that, for all $\alpha<1$, $\E |A(\ell)|^2 \lesssim |\ell|^{2\alpha}$.
Furthermore, it follows from Stokes' theorem and the fact that $\mrd A$ is a white noise, that
$A(\partial P) = \int_P \mrd A$ and hence
$\E |A(\partial P)|^2 = |P|$.
A Kolmogorov argument then implies that the GFF has a modification with $|A|_\alpha<\infty$ almost surely.
\end{remark}

To extend gauge equivalence to $\state$, let $\mfG^{0,\alpha}$ denote the closure of smooth functions in $C^{\alpha}(\T^2,G)$.
One can show that the group action $(A,g) \mapsto A^g$, defined for smooth $1$-forms and gauge transformations,
extends to a locally Lipschitz map
$
\state\times \mfG^{0,\alpha} \to \state
$,
see~\cite[Thm.~3.27, Cor.~3.36]{CCHS20}.
We then extend gauge equivalence $\sim$ to $\state$ by
\begin{equ}
A\sim B \;\; \Leftrightarrow \;\; A^g=B \text{ for some } g\in\mfG^{0,\alpha}\;.
\end{equ}
With these definitions, $\state$ has the following desirable properties (see~\cite[Thm.~2.1]{CCHS20}).

\begin{itemize}
\item For every $\ell=(x,v)\in\mcX$ and $A\in\state$, one has $|\ell_A|_{C^\alpha} \leq |\ell|^\alpha|A|_{\gr\alpha}$,
where $\ell_A\colon[0,1]\to\mfg$ is the path $\ell_A(t)=\int_0^t A(x+sv)v\mrd s$.
The holonomy $\hol(A,\ell)\in G$, defined by $\hol(A,\ell)=y_1$ where $y$ solves the ODE
\begin{equ}
\mrd y_t = y_t\mrd \ell_A(t)\;, \quad y_0=1\;,
\end{equ}
is therefore well-defined by Young integration~\cite{Young,Lyons94,FrizHairer20}.

More generally, $\hol(A,\gamma)$ is well-defined for \textit{any} $\gamma\in C^{1,\beta}([0,1], \T^2)$ where
$\beta\in (\frac2\alpha-1,1]$,
and the map $(A,\gamma)\mapsto \hol(A,\gamma)$ is H\"older continuous.
In particular classical Wilson loop observables are well-defined with good stability properties.
The relation $\sim$ can be expressed entirely in terms of holonomies.

\begin{remark}
Since $\hol(A,\gamma)$ is independent of the parametrisation of $\gamma$,
it is natural to also measure the regularity of $\gamma$ in
a parametrisation independent way.
Such a notion of regularity is introduced in~\cite[Sec.~3.2]{CCHS20} which interpolates between $C^1$ and $C^2$ (akin to how $p$-variation is a parametrisation invariant interpolation between $C^0$ and $C^1$).
\end{remark}

\item One has the embeddings
\begin{equ}
C^{\alpha/2} \hookrightarrow \state \hookrightarrow\Omega^1_{\gr\alpha} \hookrightarrow C^{\alpha-1}\;,
\end{equ}
where $\Omega^1_{\gr\alpha}$ is the completion of smooth functions under $|\cdot|_\alpha$.
(Only the last of these is non-trivial, see~\cite[Prop.~3.21]{Chevyrev19YM}.)
These embeddings are furthermore optimal in the sense that $\alpha/2$ in $C^{\alpha/2}$
cannot be decreased and $\alpha-1$ in $C^{\alpha-1}$ cannot be increased.
Remark also that $|A|_{\gr1}\asymp |A|_{L^\infty}$, while $\state$, since $\alpha<1$, contains distributions which cannot be represented by functions, such as the GFF.

\item There exists a complete metric $D$ on the quotient space of gauge orbits $\mfO = \state/{\sim}$
which induces the quotient topology.
To define $D$, we first define a new (but topologically equivalent)
metric $k$ on $\state$ by shrinking the usual metric $\Sigma(\cdot,\cdot) = |\cdot-\cdot|_\alpha$ 
in such a way that the diameter of every $R$-sphere $S_R\eqdef \{A\in\state\,:\,|A|_\alpha=R\}$
goes to zero as $R\to\infty$, but the distance between $S_r$ and $S_{R}$
for large $r\leq R$ is of order $\frac{ R}{r}-1$, so in particular goes to $\infty$ as $R\to\infty$.
Then $D$ is defined as the Hausdorff distance associated to the metric $k$ on $\state$.
\end{itemize}

\begin{remark}
The space $\state$ strengthens the definition of a Banach space $\state^\axis$ 
introduced in~\cite{Chevyrev19YM};
$\state^\axis$ is defined in a similar way
but with $\mcX$ taken as the set of \textit{axis-parallel} line segments.
The main result of~\cite{Chevyrev19YM} is that, if $G$ is simply connected, then 
there exists a (non-unique) probability measure on $\state^\axis$
such that the holonomies along all axis-parallel curves agree in distribution with those of the YM measure on $\T^2$ defined in~\cite{Sengupta97, Levy03}.
The proof of this result uses a gauge-fixed lattice approximation, which explains the restriction to axis-parallel lines.
\end{remark}


\subsection{Local solutions}\label{subsec:2D_local}

It turns out that in 2D we can sharpen the statement of Theorem~\ref{thm:local} as follows.

\begin{theorem}\label{thm:2D_local}
For every $\mathring C\in L(\mfg,\mfg)$, mollifier $\moll$, and initial condition $A(0)\in\state$,
the solution to
\begin{equ}\label{eq:SYM_2D}
\partial_t A = \Delta A +A\partial A + A^3 + (\lambda C^\eps_{\SYM} + \mathring C)A + \xi^\eps\;,
\end{equ}
converges in probability in $C(\R_+,\state\sqcup\{\skull\})$ as $\eps\downarrow0$.
The constant $C^\eps_{\SYM}$ converges to a finite limit as $\eps\downarrow0$ and is defined by
\[
C^\eps_{\SYM}=4 \hat C^\eps - \bar C^\eps\;,\quad 
\hat C^\eps = \int\partial_j K^\eps(\partial_j K *K^\eps)
\;,\quad
\bar C^\eps = \int (K^\eps)^2\;.
\]
Here $K\colon\R\times\R^2\setminus\{0\}\to\R$ is any spatially symmetric function invariant under flipping coordinates $x_i\mapsto -x_i$, $i=1,2$,
which vanishes for negative times, has bounded support, and is equal to the heat kernel $(\partial_t-\Delta)^{-1}$ in a neighbourhood of the origin.
We write
$K^\eps=\moll^\eps*K$ and $\partial_j$ is any spatial derivative, $j=1,2$.
The operator $\lambda\in L_G(\mfg,\mfg)$ is the Casimir element of $\mfg$ in the adjoint representation.

The $\eps\downarrow0$ limit of $A$, which solves~\eqref{SYM} with bare mass $\mathring C$, depends on $K$ and $\mathring C$ but not on $\moll$.
\end{theorem}

\begin{remark}
If $\mfg$ is simple (which one can assume without loss of generality, see~\cite[Remark~2.8]{CCHS20}),
then $\lambda<0$ is just a scalar.
\end{remark}

Recall from Remark~\ref{rem:2D_const_conv}
that the convergence of $C^\eps_\BPHZ \eqdef \lambda C^\eps_{\SYM}$ to a finite limit is special to dimension 2 and is due to a cancellation between the diverging constants $4\hat C^\eps$ and $\bar C^\eps$.

We briefly describe the ingredients in the
proof of Theorem~\ref{thm:2D_local},
which is based on the theory of regularity structures.
We only mention the overall strategy behind this theory, and refer to~\cite{Hairer16_CDM,ChandraWeber17,FrizHairer20} for an introduction and more details.
To solve an SPDE such as~\eqref{eq:SYM_simple},
one constructs a sufficiently large `regularity structure'
and lifts the equation to a space of `modelled distributions' with values in the regularity structure.
This construction, first introduced in~\cite{Hairer14}, is done at a high level of generality in~\cite{BHZ16}.
One then constructs a finite number of
stochastic objects from the noise called a `model' -- these objects are essentially renormalisations of functions of the form~\eqref{eq:trees_2D}-\eqref{eq:singular_funcs_2D} below,
the existence of which follows from~\cite{CH16}.
The point of the construction is that
the products $A\partial A, A^3$ and convolution with the heat kernel become stable operations on modelled distributions,
and one can solve a
fixed point problem for the `lifted' equation.
Finally, one maps the resulting modelled distribution to a distribution on $\R\times\T^2$ via the `reconstruction operator'
and identifies it with a solution to a classical renormalised PDE, at least for $\eps>0$.
This final step is carried out systematically in~\cite{BCCH21}.
All these operations are done in a way that is stable as $\eps\downarrow0$, thereby showing the desired convergence.

One of the contributions of~\cite{CCHS20}
is to develop a framework in which the algebraic results from~\cite{BHZ16,BCCH21} can be transferred to a setting in which
the noise and solution are vector-valued.
The articles~\cite{BHZ16,BCCH21} provided a general method to compute the renormalised form
of a system of \textit{scalar} SPDEs,
which in principle does apply to~\eqref{eq:SYM_simple} by writing it as a system of $d\times \dim(\mfg)$ scalar-valued equations using a basis.
However, such a procedure is cumbersome and	unnatural; it is more desirable to find 
a framework that preserves the vector-valued nature of the noise and solution, which is the purpose of~\cite[Sec.~5]{CCHS20}.

The main idea behind the extension in~\cite{CCHS20} is to define a category of `symmetric sets' and a functor between this category and the category of vector spaces.
This construction allows one to canonically associate partially symmetrised tensor products of vector spaces to combinatorial rooted trees which commonly appear in regularity structures.
One of the main outcomes is a procedure to compute the renormalised form of equations like~\eqref{eq:SYM_simple} without resorting to a basis.

In working out the renormalised equation~\eqref{eq:SYM_2D}, it is
easy to deduce by power-counting and symmetry arguments that the only non-vanishing counterterms arise from the three trees
%
\begin{equ}\label{eq:trees_2D}
\<I[I'Xi]I'Xi_notriangle>\;,\quad
\<IXiI'[I'Xi]_notriangle>\;,\quad\textnormal{and}\quad 
\<IXi^2>\;.
\end{equ}
Before renormalisation, these trees correspond respectively to the three functions
\begin{equ}\label{eq:singular_funcs_2D}
[K*(\partial_j K*\xi)] \cdot (\partial_j K*\xi)\;, 
\quad
 (K*\xi) \cdot [\partial_j K*(\partial_j K*\xi)]\;,
\quad\textnormal{and}\quad
(K*\xi)^2\;.
\end{equ}
A computation shows that the final counterterm is precisely $\lambda C^\eps_\SYM A$ as defined in Theorem~\ref {thm:2D_local} (see~\cite[Lem.~6.14]{CCHS20}).

It follows from the general theory of regularity structures that~\eqref{eq:SYM_2D}, for any $\mathring C\in L(\mfg,\mfg)$, converges locally in time in $C^{\alpha-1}$.
To improve this to convergence in $C(\R_+,\state\sqcup\{\skull\})$, one decomposes the solution $A$ into $A=\Psi+B$
where $\Psi$ solves the SHE $\partial_t \Psi= \Delta\Psi + \xi$ with initial condition $\Psi(0)=A(0)$
and $B$ is in $C^{1-\kappa}$ for ay $\kappa>0$.
One can then show by hand that $\Psi\in C(\R_+,\state)$ (see~\cite[Sec.~4]{CCHS20}),
which, together with the embeddings $C^{\alpha/2}\hookrightarrow \state\hookrightarrow C^{\alpha-1}$,
shows that $A$ converges to a maximal solution with values in $\state$.

\subsection{Gauge covariance}
\label{subsec:covar_2D}

Recall that Theorems~\ref{thm:generative} and~\ref{thm:B_bar_A} imply a form of gauge covariance for~\eqref{SYM}.
Theorem~\ref{thm:generative}\ref{pt:gen_exists}
is a relatively straightforward consequence of Theorem~\ref{thm:2D_local}.
One defines the random variable $A$ by solving~\eqref{SYM} until the first time that $|A(t)|_\alpha \geq 2+\inf_{B\sim A(t)} |B|_\alpha$,
at which point one uses a measurable selection $S\colon\mfO\to\state$ to jump to a new small representative $B$ of the gauge orbit $[A(t)]$ for which $|B|_\alpha < 1+ \inf_{a\in [A(t)]} |a|_\alpha$.
These jump times define the increasing sequence of stopping times $\{\varsigma_j\}_{j\geq 0}$ in item~\ref{pt:SYM_jumps}.
Items~\ref{pt:init_cond}-\ref{pt:honest_blow_up} all follow readily from the construction.

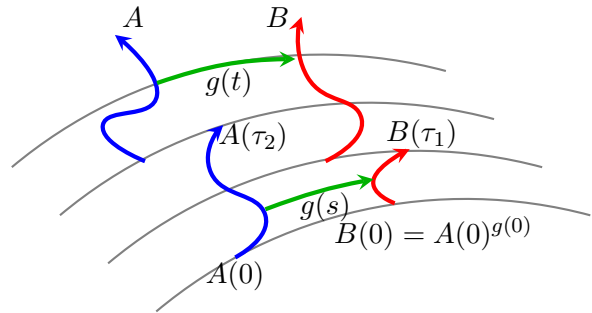
\begin{figure}[h]
\centering
  \begin{tikzpicture}[scale=0.8]
  
    \draw[thick,gray] (0,2.4) arc (75:130:8cm);
    \draw[thick,gray] (.8,1.6) arc (75:130:8cm);
    \draw[thick,gray] (1.6,0.8) arc (75:130:8cm);
    \draw[thick,gray] (2.4,0) arc (75:130:8cm);
    
 \draw [blue,ultra thick,->] plot [smooth, tension=1] coordinates 
 { (-3.5,-0.7) (-3.0,0) (-3.9, 0.7) (-3.7,1.5)};
\draw [blue,ultra thick,->] plot [smooth, tension=1] coordinates 
 { (-5.0,0.9) (-5.7 , 1.5) (-4.8 , 2) (-5.5,  3)};   
\node at (-3.5, -1) {$A(0)$}; 
\node at (-5.2,  3.3) {$A$}; 
\node at (-3.2,1.3) {$A(\tau_2)$}; 
     
 \draw [red,ultra thick,->] plot [smooth, tension=1] coordinates 
 { (-0.85,0.2) (-1.2,0.6) (-0.6,1.1)};
\draw [red,ultra thick,->] plot [smooth, tension=1] coordinates 
 { (-2,0.9) (-1.4,1.6) (-2.3 , 2.3) (-2.4,3.3)};
\node at (-0.2, -0.3) {$B(0)= A(0)^{g(0)}$}; 
\node at (-2.8,  3.3) {$B$}; 
\node at (-0.4,1.35) {$B(\tau_1)$};
 
   

\draw [darkgreen, ultra thick,-> ,bend left =5]  (-3,0.1)  to (-1.2,0.6);
 \node at (-2,0.15) {$g(s)$}; 
 
\draw [darkgreen, ultra thick,-> ,bend left =9]  (-4.8,  2.2)  to (-2.5, 2.6);
 \node at (-3.6,2.2) {$g(t)$};  

\end{tikzpicture}
\caption{Two random variables $A,B$ associated to generative probability measures $\mu,\nu$ respectively.
The set of jump times $\{\tau_1,\tau_2,\ldots\}$
is the union of the sets of stopping times $\{\varsigma_1,\varsigma_2,\ldots\}$ and $\{\sigma_1,\sigma_2,\ldots\}$ associated to $\mu,\nu$.
Grey lines indicate gauge orbits.
There is a coupling such that $A^g=B$ for a time-dependent gauge transformation $g$.
In between the jump times, $A$ and $B$ solve~\eqref{SYM} with bare mass $\check C$ driven by $\xi$ and $\bar\xi=\Ad_g\xi$ respectively.
}
\label{fig:Markov}
\end{figure}

The proof of Theorem~\ref{thm:generative}\ref{pt:covar}, which is the main statement of Theorem~\ref{thm:generative},
requires more work.
The idea is to use Theorem~\ref{thm:B_bar_A},
which we admit for now, to couple the solutions to~\eqref{SYM}
with bare mass $\check C$ and initial conditions $a\sim b$.
Specifically, let $\nu$ be a generative probability measure with bare mass $\check C$ and initial condition $b\in\state$, and consider any $a\sim b$.
Letting $B$ and $\bar \xi$ denote the random variable and white noise respectively corresponding to $\nu$,
it follows from Theorem~\ref{thm:B_bar_A} that, on the same probability space,
there exists a c\`adl\`ag process $A$ defined as above with~\eqref{SYM} driven by $\xi\eqdef \Ad_{g^{-1}}\bar \xi$ and bare mass $\check C$ in such a way that $B=A^g$.
Here $g$ is c{\`a}dl{\`a}g with values in $\mfG^{0,\alpha}\sqcup\{\skull\}$ and jump times contained in those of $A$ and $B$, and solves~\eqref{eq:SPDE_B} in between these jump times, see Figure~\ref{fig:Markov}.
This shows that the pushforward of $\nu$ to the orbit space $\mfO$ is equal to the 
pushforward of the law of $A$ from the proof of Theorem~\ref{thm:generative}\ref{pt:gen_exists},
which proves  Theorem~\ref{thm:generative}\ref{pt:covar}.

\begin{remark}
A crucial fact used in the above construction is that the solution $g$ to~\eqref{eq:SPDE_B} does not blow up before $B$ or $A\eqdef B^{g^{-1}}$ does.
This is due to the elementary but important property of $\state$ that
\begin{equ}\label{eq:g_bound_2D}
|g|_{C^\alpha} \lesssim 1+|A|_{\gr\alpha} + |A^g|_{\gr\alpha}\;,
\end{equ}
which in turn follows from the facts that, for any curve $\gamma_{xy}$ with $\gamma_{xy}(0)=x$ to $\gamma_{xy}(1)=y$,
\begin{equ}\label{eq:g_identity_hol}
g(y)=\hol(A^g,\gamma_{xy})^{-1}g(x)\hol(A,\gamma_{xy})\;,
\end{equ}
and that, if $\gamma_{xy}$ is the shortest such curve, then the distance of the holonomy $\hol(A,\gamma_{xy})\in G$ from the identity in $G$
is of order $|x-y|^\alpha|A|_{\gr\alpha}$ by Young ODE theory (see~\cite[Sec.~3.5]{CCHS20}).
\end{remark}

The proof of Theorem~\ref{thm:generative}\ref{pt:covar} therefore boils down to that of Theorem~\ref{thm:B_bar_A},
which is the main content of~\cite[Sec.~7]{CCHS20}.
The proof of Theorem~\ref{thm:B_bar_A} proceeds by taking a non-anticipative mollifier $\moll$ and considering the two systems
\begin{equs}[eq:SPDE_B_bares]
\partial_t B &= \Delta B + B\partial B + B^3
+\Ad_g \xi^\eps
+ \mathring C_1 B + \mathring C_2 (\mrd g) g^{-1}\;,
\\
(\partial_t g)g^{-1} &= \partial_j((\partial_j g )g^{-1}) + [B_j,(\partial_j g) g^{-1}]\;,
\end{equs}
and
\begin{equs}[eq:SPDE_bar_A_bares]
\partial_t \bar A &= \Delta \bar A + \bar A\partial \bar A + \bar A^3
+\moll^\eps*(\Ad_{\bar g} \xi)
+ \mathring C_1 \bar A + \mathring C_2 (\mrd \bar g) \bar g^{-1}\;,
\\
(\partial_t \bar g)\bar g^{-1} &= \partial_j((\partial_j \bar g )\bar g^{-1}) + [\bar A_j,(\partial_j \bar g) \bar g^{-1}]\;,
\end{equs}
where $\mathring C_1,\mathring C_2\in L(\mfg,\mfg)$ are new arbitrary bare masses.
These two systems represent the unrenormalised versions of~\eqref{eq:SPDE_B} and~\eqref{eq:SPDE_bar_A} respectively with an extra free parameter (the bare mass $\mathring C_2$).
(In~\cite{CCHS20} one actually writes the above systems
using new variables $U= \Ad_g$ and $h=(\mrd g)g^{-1}$,
which assists in computing the desired renormalised equations.)

\begin{remark}\label{rem:input_dist_2D}
Well-posedness for the systems~\eqref{eq:SPDE_B_bares}-\eqref{eq:SPDE_bar_A_bares} as $\eps\downarrow0$
is generally standard.
However, a subtlety arises from the multiplicative noise term $\Ad_g\xi$ which
is in $C^{-2-\kappa}$ and leads to problems in posing a suitable fixed point problem
($-2$ is the threshold regularity at which one cannot extend uniquely a distribution from $\R_+\times\R^{d}$ to $\R\times\R^d$).
This issue is handled by decomposing
$g=\CP g(0) + \hat g$, where $\CP g(0)$ is the harmonic extension of $g(0)$ to positive times.
Then $\hat g$ vanishes at $t=0$, which makes the product $\Ad_{\hat g}\xi$ better behaved,
while the product $\Ad_{\CP g(0)} \xi$
is shown to be a well-defined distribution in $C^{-2-\kappa}(\R\times\R^d)$ using stochastic estimates.
\end{remark}

To show that some renormalised forms of~\eqref{eq:SPDE_B_bares} and~\eqref{eq:SPDE_bar_A_bares} converge to the \textit{same} limit as $\eps\downarrow0$,
the strategy taken in~\cite{CCHS20}
is to introduce $\eps$-dependent norms on the underlying regularity structure.
The idea behind these norms is that they allow one to lift the basic estimate
\[
|f - \moll^\eps*f|_{C^\ell} \lesssim \eps^{\theta}|f|_{C^{\ell+\theta}}
\]
to the level of modelled distributions.
It is then possible to show that the modelled distributions which solve the systems~\eqref{eq:SPDE_B_bares} and~\eqref{eq:SPDE_bar_A_bares}
are, for small times, at distance of order $\eps^\theta$ for some $\theta>0$ sufficiently small.
By continuity of the reconstruction operator,
it follows that some renormalised forms of~\eqref{eq:SPDE_B_bares} and~\eqref{eq:SPDE_bar_A_bares} converge to the same limit as $\eps\downarrow0$.


We next identify the renormalised forms of~\eqref{eq:SPDE_B_bares} and~\eqref{eq:SPDE_bar_A_bares}.
We drop for the moment the assumption that $\moll$ is non-anticipative since this step can be done without it.
It follows from a direct computation using the algebraic theory developed in~\cite{BCCH21} and~\cite[Sec.~5]{CCHS20}
that the renormalised equations are
\begin{equ}\label{eq:B_bare_mases_renorm}
\partial_t B = \Delta B + B\partial B + B^3  + ( \lambda C^\eps_{\SYM} +\mathring C_1) B + \Ad_g(\xi^\eps) + (\lambda\tilde C^\eps + \mathring C_2) (\mrd g)g^{-1}
\end{equ}
and
\begin{equ}\label{eq:bar_A_bare_mases_renorm}
\partial_t \bar A = -\Delta \bar A + \bar A\partial \bar A + \bar A^3 
+ (\lambda C^\eps_{\SYM} +\mathring C_1) \bar A
+ \moll^\eps*(\Ad_g \xi) + (\lambda\tilde C^{0,\eps} + \mathring C_2) (\mrd \bar g)\bar g^{-1}
\end{equ}
where 
\begin{equ}
\tilde C^\eps = \int \moll^\eps (K*K^\eps)\;,
\quad \tilde
C^{0,\eps}=\lim_{\delta\downarrow0} \int \moll^\delta (K * \moll^\delta*K^\eps) = (K*K^\eps)(0)\;,
\end{equ}
and $C^\eps_{\SYM}$ is the same constant appearing in Theorem~\ref{thm:2D_local}.
Remark that $\tilde C^\eps$ and $\tilde C^{0,\eps}$ both converge to finite limits as $\eps\downarrow0$ because, if $G$ is the heat kernel,
then $(G*G)(t,\cdot)=tG(t,\cdot)$, which is a bounded functions because $d=2$.
The components $g$ and $\bar g$ do not require renormalisation and solve the same equations as in~\eqref{eq:SPDE_B_bares} and~\eqref{eq:SPDE_bar_A_bares}.

The values of $\tilde C^\eps$ and $\tilde C^{0,\eps}$
are respectively the BPHZ constants associated with the two trees
\begin{equ}\label{eq:2D_finite_trees}
\<I[I[Xi]]Xi_typed>
\quad \text{and} \quad
\<I[I[Xi]]Xi_typed*>\;,
\end{equ}
which correspond, before renormalisation,
to the two functions
\begin{equ}
\xi^\eps \cdot (K*K*\xi^\eps) 
\quad \text{and} \quad
\xi \cdot (K*K^\eps*\xi)\;.
\end{equ}

\begin{remark}
To derive and solve the equation for $\bar A$, one substitutes $\xi$ by $\xi^\delta$
and takes the limit $\delta\downarrow0$ with $\eps>0$ fixed;
this ensures that all objects are smooth for $\eps,\delta>0$.
This also explains the definition of $\tilde C^{0,\eps}$ as a limit in $\delta\downarrow0$.
\end{remark}

To summarise, for any $\mathring C_1,\mathring C_2\in L(\mfg,\mfg)$ and any mollifier $\moll$,
the solutions $B$ and $\bar A$ to~\eqref{eq:B_bare_mases_renorm} and~\eqref{eq:bar_A_bare_mases_renorm} respectively
converge to the same limit as $\eps\downarrow0$ over a short random time interval (the argument in~\cite{CCHS20} that $(B,g)$ and $(\bar A,\bar g)$ converge as \textit{maximal} solutions uses non-anticipativity of $\moll$).

To make the identification with the equations in Theorem~\ref{thm:B_bar_A},
we take $\mathring C_1= \mathring C$.
We then set $\mathring C_2=\lim_{\eps\downarrow0}\lambda C^\eps_\SYM + \mathring C-\lambda \tilde C^\eps$,
so that
$\lambda \tilde C^\eps +\mathring C_2=\lambda C^\eps_\SYM + \mathring C + o(1)$ as $\eps\downarrow0$.
Finally, we want to find $C$ (playing the role of $\check C$) so that
$\lambda \tilde C^{0,\eps} + \mathring C_2 = \mathring C_1 - C + o(1)$.
This desired operator is
\begin{equ}\label{eq:C_def_2D}
C
= \lim_{\eps \downarrow 0} \lambda(\tilde C^\eps-C^\eps_\SYM - \tilde C^{0,\eps})\;.
\end{equ}
We have now argued that the solutions $B$ and $\bar A$ as in Theorem~\ref{thm:B_bar_A}
with $\check C=C$
converge to the same limit over a short time interval as $\eps\downarrow0$ for \textit{any} mollifier $\moll$.

Next, we claim that if $\moll$ is \textit{non-anticipative},
then $C$ is independent of $\moll$.
Indeed, $\tilde C^{0,\eps}=0$ for non-anticipative $\moll$.
Furthermore, it follows from the identity $(\partial_t-\Delta)K = \delta + Q$, where $\delta$ is the Dirac delta and $Q$ is smooth and supported away from the origin,
and a computation with integration by parts (see~\cite[Lem.~6.9]{CCHS20}), that $C$ in this case is equal to
\begin{equ}\label{eq:check_C_2D}
\check C = \lim_{\eps \downarrow 0} \lambda(\tilde C^\eps-C^\eps_\SYM)
=\lim_{\eps\downarrow0}\lambda
\Big(-\int (K*K^\eps)(Q*\moll^\eps) - \int (Q*K^\eps) K^\eps\Big)\;.
\end{equ}
Since $Q$ is supported away from the origin,
the final limit is independent of $\moll$.

To conclude the proof of Theorem~\ref{thm:B_bar_A},
it remains to show that~\eqref{SYM} with bare mass $\check C$ defined by~\eqref{eq:check_C_2D}
is independent of $\moll$ and of $K$.
Remark that $\moll$ can now be \textit{any} mollifier, not necessarily non-anticipative.
Independence of $\moll$ follows from the final part of Theorem~\ref{thm:2D_local} since $\check C$ is independent of $\moll$.
Independence of $K$ follows from the fact that $\lim_{\eps\downarrow0}(\lambda C^\eps_\SYM + \check C)=\lim_{\eps\downarrow0}\lambda \tilde C^\eps$, which does not depend on the choice of $K$
since $K$ is always equal to the heat kernel near the origin.

\begin{remark}
The final mass renormalisation one takes in~\eqref{SYM} to obtain gauge covariance in 2D is therefore
\begin{equ}\label{eq:bar_C_2D}
\lim_{\eps\downarrow0}(\lambda C^\eps_{\SYM} + \check C) = \lim_{\eps\downarrow0} \lambda\tilde C^\eps = \lim_{\eps\downarrow0}  \lambda \int \moll^\eps (K*\moll^\eps*K)\;,
\end{equ}
which is the constant `$\bar C$' appearing in~\cite[Thm.~2.9]{CCHS20}.
\end{remark}

\begin{remark}
The existence of $\check C$ with the above properties may appear as a bit of a miracle.
Indeed, the fact that $\tilde C^\eps$ and $\tilde C^{0,\eps}$ converge to finite limits was easy to see
because $(G*G)(t,\cdot)=tG(t,\cdot)$, is a bounded function for the heat kernel $G$ in 2D.
On the other hand,
the fact that $C^\eps_\SYM$ converges to a finite limit, and thus that $\check C$ exists, is not a priori obvious because it relies on a cancellation between diverging BPHZ constants $4\hat C^\eps$ and $\bar C^\eps$ in Theorem~\ref{thm:2D_local}.
The fact that $\check C$ is furthermore independent of $\moll$ relies on a cancellation between $\tilde C^\eps$ and $C^\eps_\SYM$, and that~\eqref{SYM} with bare mass $\check C$ is independent of $K$ relies on the expression for $\tilde C^\eps$. 


These cancellations, convergences, and independencies are shown in~\cite{CCHS20}
using explicit computations.
We will see in Section~\ref{subsec:3D_gauge_covar}
that there is another argument which allows us to see the existence of $\check C$ by instead exploiting symmetries of the equation.
\end{remark}

\section{Three dimensions}\label{sec:3D}

We now discuss the main results of~\cite{CCHS22}, which deals with the 3D theory.

\subsection{State space}\label{subsec:3D_state}

The idea of defining a state space of $1$-forms in terms of line integrals no longer works in 3D because already the GFF $A$ is too singular to be restricted to lines.
To see this, recall that the correlation function of $A$ in 3D behaves like $C(x,y) \asymp \frac{1}{|x-y|}$ for $x,y$ close (vs. $C(x,y) \asymp -\log|x-y|$ in 2D).
Therefore, for $\ell=(x,v)$,
\begin{equ}
\E|A(\ell)|^2
= \int_0^1\mrd t\int_0^1\mrd s |\ell|^2C(x+tv, x+sv)
\asymp \int_0^1\mrd t\int_0^1\mrd s \frac{|\ell|}{|t-s|} = \infty\;.
\end{equ}
This strongly suggests it is hopeless to find a state space for 3D quantum Yang-Mills such that 
holonomies are well-defined as we did in 2D.

The construction of the state space $\state$ in~\cite{CCHS22} proceeds in two steps.
The first step is to define a space $\init$ of initial conditions for a gauge-covariant regularising operator.
Abstractly, we will find a metric space of distributional $1$-forms $(\init,\Theta)$
and a family of operators $\{\mcF_t\}_{t > 0}$, $\mcF_t\colon \mcI\to C^\infty$, such that
\begin{enumerate}[label=(\alph*)]
\item \label{pt:F_gauge_covar}for smooth $A,B$, 
\begin{equ}\label{eq:sim1}
A\sim B \;\;\Leftrightarrow \;\;
\mcF_t(A) \sim \mcF_t(B) \;\; \text{for some } t>0\;,
\end{equ}

\item \label{pt:F_cont} $\mcF_t \colon \mcI\to C^\infty$ is continuous for every $t>0$.
\end{enumerate}
If we can find such $\init$ and $\{\mcF_t\}_{t>0}$, then we can extend gauge equivalence $\sim$ to $\init$ by using~\eqref{eq:sim1} as a \textit{definition}.
Finally, we want $\init$ to be sufficiently large to contain distributions as rough as the GFF.

Our concrete choice for $\mcF_t(a)$ is the solution $A(t)$ to the deterministic YM flow (with DeTurck term)\footnote{To avoid arguing that~\eqref{eq:YM_flow} has global in time solutions, the definition of $\mcF_t$ in~\cite{CCHS22} is restricted to short intervals $t\in(0,T)$ -- see Proposition~\ref{prop:YM_extends} below --
and items~\ref{pt:F_gauge_covar} and~\ref{pt:F_cont} should also be understood locally in time.
To simplify the exposition, we ignore this detail here.}
\begin{equ}\label{eq:YM_flow}
\partial_t A = -\mrd^*_A F(A) - \mrd_A \mrd^* A = \Delta A + A\partial A+ A^3\;,\quad A(0)=a\;.
\end{equ}
This choice is natural because we ultimately want to start~\eqref{SYM}, the stochastic version of~\eqref{eq:YM_flow},
from initial data as rough as the GFF, which is at least as hard as solving~\eqref{eq:YM_flow}.

\begin{remark}
The fact that~\ref{pt:F_gauge_covar} above holds for $\mcF$ was likely already known in the literature, but a proof is given in~\cite[Sec.~2.2]{CCHS22} based on analytic continuation
and global existence of the YM flow without DeTurck term~\cite{Rade92,HT04}.
\end{remark}

The second step in the construction of $\state$ is to augment $\init$ with an additional norm which ensures that a form of the bound~\eqref{eq:g_bound_2D} holds.
This turns out to be critical in several places of the construction for the Markov process in Theorem~\ref{thm:generative}\ref{pt:covar}.

We mention that the idea to use the YM flow to define a suitable space of distributional $1$-forms
was already suggested in~\cite{CG13} (see also~\cite{NN06,Luscher10,Fodor12, MR3133916} for related ideas in physics).
We also point out that another state space which bears close similarity to $\init$
was independently defined in~\cite{Sourav_state,Sourav_flow}
and was shown to support the GFF.

\subsubsection{The first half}
\label{subsubsec:1st_half}

Before defining $\init$, we briefly illustrate why classical H\"older--Besov spaces $C^\eta$ are not suitable for our purposes.
The standard strategy to solve~\eqref{eq:YM_flow} is to rewrite the equation in mild formulation
\begin{equ}
A(t) = \CM_t(A)\eqdef \CP_t a + \int_0^t \CP_{t-s}[A(s)\partial A(s) + A(s)^3]\mrd s\;,
\end{equ}
where $\CP_t=e^{t\Delta}$ is the heat flow, and
show that the map $\CM\colon A \mapsto \CM(A)$
is a contraction on a ball in a suitable Banach space.
This Banach space should at least
contain $\CM(0) = \CP a$ and $\CM(\CP a)$ restricted to a short time interval,
which are the first and second Picard iterates respectively.
However, for generic $a\in C^\eta$, the best we can do to handle the product $\CP_{t-s}[\CP_s a\cdot \partial \CP_s a]$ in $\CM(\CP a)$ is to
estimate for $s>0$
\begin{equs}{}
&|\CP_s a|_{L^\infty}
\lesssim s^{\frac\eta2}|a|_{C^\eta}
\;,\quad
|\partial \CP_s a|_{L^\infty} \lesssim s^{\frac\eta2 -\frac12}|a|_{C^\eta} \label{eq:a_Shauder}
\\
&\Rightarrow
|\CP_s a \cdot \partial \CP_s a|_{L^\infty} \lesssim s^{\eta-\frac12}|a|_{C^\eta}\;. \label{eq:a_prod}
\end{equs}
We now recall that the GFF takes values in $C^\eta$ for any $\eta<-\frac12$
(and not for $\eta=-\frac12$).
Therefore, trying to estimate the term $\int_0^t \CP_{t-s}[\CP_s a \cdot \partial \CP_s a]\mrd s$ in the 2nd Picard iterate leads to a non-integrable singularity $\int_0^t s^{\eta-\frac12} \mrd s = \infty$.
The above estimates are, in general, sharp, which suggests that one cannot start the YM flow~\eqref{eq:YM_flow} from initial data in $C^\eta$ with $\eta<-\frac12$ (and even in $C^{-1/2}$).

To circumvent this issue and motivate the definition of $\init$, remark that the GFF $a$ is a highly non-generic element of $C^\eta$
and  $\CP_s a \cdot \partial \CP_s a$ behaves better than the naive bound~\eqref{eq:a_prod} would suggest.
Indeed, a 2nd moment estimate shows that for any $\beta\in(-1,0)$ and $\delta>1+\frac\beta2$, almost surely uniformly in $s\in(0,1)$
\begin{equ}\label{eq:GFF_prod}
|\CP_s a \cdot \partial \CP_s a|_{C^\beta} \lesssim s^{-\delta}\;.
\end{equ}
(The way to guess the bound~\eqref{eq:GFF_prod} is to use the estimates $|\CP_s a|_{C^{\beta/2}}
\lesssim s^{-\frac14-\kappa-\frac\beta4}$
and
$|\partial \CP_s a|_{C^{\beta/2}} \lesssim s^{-\frac34-\kappa - \frac\beta4}$ for any $\kappa>0$,
and then pretend that multiplication is a bounded operator $C^{\beta/2}\times C^{\beta/2}\to C^{\beta}$ -- this is clearly false since $\beta<0$, but ends up working for $(\CP_s a,\partial \CP_s a)$ due to probabilistic cancellations.)

Since we can take $\delta<1$ in~\eqref{eq:GFF_prod} (vs. $-\eta+\frac12>1$ in~\eqref{eq:a_prod}), this improved regularity ends up being enough to show that every Picard iterate of $\CM$ is well-defined when $a$ is the GFF.
It is therefore natural to make the following definition.

\begin{definition}
For $\eta<-\frac12$, $\beta<0$, and $\delta\in(1+\frac\beta2,1)$,
let $\init$ be the completion of smooth $1$-forms under the metric
\[
\Theta(A,B) \eqdef |A-B|_{C^\eta} + \sup_{t\in(0,1)} t^\delta|\CP_t A\cdot\partial \CP_t A - \CP_tB\cdot\partial \CP_tB|_{C^{\beta}}\;.
\]
\end{definition}

\begin{remark}
The space $\init$ can be identified with a subset of $C^{0,\eta}\eqdef \overline{C^\infty}^{C^\eta}$
because the map $C^\eta\ni A\mapsto \CP_\cdot A \cdot \partial \CP_\cdot A \in C^\infty((0,1),C^\beta)$ has a closed graph.
\end{remark}

A standard argument with Young's product theorem and estimates of the type~\eqref{eq:a_Shauder} shows that the YM flow extends to $\init$ in the following sense.

\begin{proposition}\label{prop:YM_extends}
For every ball $\mcB$ in $(\init,\Theta)$ centred at $0$, there exists $T>0$ such that for all $t\in(0,T)$, the YM flow~\eqref{eq:YM_flow} extends to a continuous function $\mcF_t\colon \mcB\to C^\infty$ (which is Lipschitz for any norm on $C^\infty$).
\end{proposition}
We can therefore extend gauge equivalence $\sim$ to $\init$ by using~\eqref{eq:sim1} as a definition.
We emphasise that $\init$ is not a vector space and that this is unavoidable:
there is no Banach space that carries the GFF on $\T^3$ and to which the YM flow~\eqref{eq:YM_flow} extends as a continuous function (locally in time), see~\cite{Chevyrev22_norm_inf}.


\subsubsection{The second half}\label{subsubsec:2nd_half}

We now refine the space $\mcI$ in order to obtain control on gauge transformations of the type~\eqref{eq:g_bound_2D},
which proves crucial in the construction of the Markov process on gauge orbits associated to~\eqref{SYM}.

\begin{definition}
For $\alpha\in(0,\frac12)$ and $\theta>0$,
define the norm $\heatgr{\cdot}_{\alpha,\theta}$ on smooth $1$-forms
\[
\heatgr{A}_{\alpha,\theta} \eqdef \sup_{t\in(0,1)} \sup_{|\ell|<t^\theta} \frac{|(\CP_t A)(\ell)|}{|\ell|^\alpha}\;,
\]
where the second supremum is over all line segments $\ell=(x,v)\in\T^3\times \{v\in\R^3\,:\,|v|<\frac14\}$ with $|\ell|\eqdef|v|<t^\theta$,
and $(\CP_tA)(\ell)$ is the line integral of $\CP_t A$ along $\ell$ defined analogously to~\eqref{eq:line_int_def}.
Define further the metric
\[
\Sigma(A,B) = \Theta(A,B) + \heatgr{A-B}_{\alpha,\theta}\;,
\]
and let $\state$ by the completion of smooth $1$-forms under $\Sigma$.
\end{definition}

To motivate the norm $\heatgr{\cdot}_{\alpha,\theta}$, recall that the estimate~\eqref{eq:g_bound_2D} relies on the identity~\eqref{eq:g_identity_hol},
which in turn requires that line integrals and holonomies of $A$ are well-defined.
But we saw that the GFF $A$ cannot even be restricted to lines!

The idea above is to
consider, instead of $A$, the heat flow regularisation $\{\CP_t A\}_{t\in(0,1)}$.
A quick computation shows that, uniformly in $0<\sqrt t< |\ell|<\frac14$,
\begin{equ}
\E |(\CP_tA)(\ell)|^2\asymp |\ell|\log(|\ell| t^{-1})
\end{equ}
(with even better bounds for $|\ell|\leq \sqrt t$).
Therefore, as expected, $\E|(\CP_tA)(\ell)|^2$ blows up for fixed $\ell$ as $t\downarrow0$, but rather slowly.
Furthermore, restricting to short length scales, say $|\ell|<t^\theta$ for any $\theta>0$,
\begin{equ}
\E |(\CP_tA)(\ell)|^2\lesssim -|\ell|\log |\ell|
\end{equ}
uniformly in $t\in(0,1)$ and $|\ell|<t^\theta$.
Combined with a Kolmogorov argument,
this shows that $\heatgr{A}_{\alpha,\theta}<\infty$ almost surely.
The restriction to $\alpha<\frac12$ is natural because the GFF in 3D has $\frac12$ less regularity than in 2D (e.g.\ $C^{-\kappa-1/2}$ in 3D vs. $C^{-\kappa}$ in 2D for H{\"o}lder--Besov regularity),
and we saw that $|A|_{\gr\alpha}<\infty$ for $\alpha<1$ for the GFF $A$ in 2D.

\begin{remark}
The metric space $\state$ can be identified with a subset of $\init\subset C^{0,\eta}$
and comes with the parameters $(\eta,\beta,\delta,\alpha,\theta)$, the possible range of which is given in~\cite[Sec.~5]{CCHS22};
$(\state,\Sigma)$ is the space appearing in Theorem~\ref{thm:state_space} for $d=3$.
\end{remark}

With these definitions, one can show that the group of sufficiently smooth gauge transformations acts continuously on $\state$.
Namely, there exists $\rho\in(\frac12,1)$ such that $(A,g) \mapsto A^g$ extends continuously to a map $\state\times \mfG^{\rho} \to \state$
where
\begin{equ}
\mfG^\rho \eqdef  C^{\rho}(\T^3,G)\;.
\end{equ}
Furthermore, this action preserves $\sim$ defined by~\eqref{eq:sim1}, i.e.\ $A\sim A^g$ for all $(A,g)\in\state\times \mfG^\rho$.

We can now state~\cite[Thm.~2.39]{CCHS22}, which is one of the main results of~\cite[Sec.~2]{CCHS22}
and the motivation behind the norm $\heatgr{\cdot}_{\alpha,\theta}$.

\begin{theorem}\label{thm:g_bound}
There exist constants $C,q >0$ and $\nu\in(0,\frac12)$ such that, for all $g\in \mfG^\rho$ and $A\in\state$,
\begin{equ}
|g|_{C^\nu} \leq C(1+\Sigma(A,0) + \Sigma(A^g,0))^q\;.
\end{equ}
\end{theorem}

The proof of Theorem~\ref{thm:g_bound} relies on two estimates:
(i) the estimate~\eqref{eq:g_bound_2D} used in the 2D case (which of course holds in arbitrary dimension),
and (ii) a `backwards estimate' which controls the initial condition of a parabolic PDE in
terms of its behaviour for positive times (see~\cite[Lem.~2.46(b)]{CCHS22}) -- this estimate is applied to the harmonic map flow-type PDE solved by $h$ for which $h(0)=g$ and $\mcF_t(A)^{h(t)} = \mcF_t(A^g)$ for all $t>0$.
Theorem~\ref{thm:g_bound} is then obtained by suitably interpolating between estimates (i) and (ii).

\begin{remark}
One can show that mollifications of the SHE converge in probability in the space $C(\R_+,\state)$  (see~\cite[Cor.~3.14]{CCHS22}).
In particular, the SHE admits a modification with sample paths in $C(\R_+,\state)$.
\end{remark}

\begin{remark}\label{rem:no_gauge_group}
Unlike in 2D, the action of $\mfG^\rho$ on $\state$ is not transitive over the orbits and it is unclear if $\sim$, or some variant of it,
is determined by the action of a group.
This lack of a gauge group is responsible for the gap in our understanding of the quotient space $\state/{\sim}$ in 3D vs. 2D, see Remark~\ref{rem:quotient_top}.
\end{remark}

\subsection{Local solutions}
\label{subsec:local_sol_3D}

We next explain how one proves Theorem~\ref{thm:local} in 3D, which is done in~\cite[Sec.~5]{CCHS22}.
We do not restate the result here like we did in Section~\ref{subsec:2D_local}
since we can't make it substantially more precise.

Though primarily using the theory of regularity structures as before,
there are two main additional challenges on top of the 2D case.
The first is purely algebraic and concerns
showing that the renormalisation counterterms are of the form $C^\eps_\BPHZ A$.
The difficulty is that there are dozens of trees which potentially contribute to renormalisation (vs. just 9 trees in 2D, see~\cite[Sec.~6.2.3]{CCHS20}).

By power counting, one can deduce that the renormalisation is \textit{linear} in $A$.
To argue that one sees precisely $C^\eps_\BPHZ\in L_G(\mfg,\mfg)$ requires a systematic approach to symmetry arguments, which is developed in~\cite[Sec.~4]{CCHS22}
and which could of independent interest in other contexts.

To give an example of how this works, we argue that the renormalisation is `block diagonal',
i.e.\ if the counterterm $c A_j$ appears in the $A_i$ equation for $j\neq i$ then $c=0$.
Indeed, if we flip the coordinate $x_i\mapsto -x_i$ and thus $\partial_i\mapsto -\partial_i$, together with $A_i\mapsto -A_i$ and $\xi_i^\eps\mapsto -\xi_i^\eps$, while keeping all terms with indexes $j\neq i$ the same,
it is immediate that all the terms in the $A_i$ equation~\eqref{eq:SYM_coord} flip sign.
Using the symmetry of the noise $\xi_i^\eps\eqlaw -\xi_i^\eps$ and invariance under the flip $x_i\mapsto -x_i$ of the kernel $K$ used to define $C^\eps_\BPHZ$,
one can show that renormalised equation must possess the same symmetry, namely all
terms in the renormalised $A_i$ equation must flip sign.
Since we kept $A_j$ the same, this shows that any factor $cA_j$ in the renormalised $A_i$ equation must have $c=0$.

The way one argues that the \textit{same} $C^\eps_\BPHZ$ appears for all $i\in\{1,2,3\}$
and that $C^\eps_\BPHZ \in L_G(\mfg,\mfg)$
is similar: one exploits symmetry under reflections $x_i\leftrightarrow x_j$ and $\xi_i^\eps\eqlaw \xi_j^\eps$ for the former, and
symmetry under constant gauge transformations $A_i \mapsto \Ad_g A_i$ and $\xi_i\eqlaw \Ad_g\xi_i^\eps$ where $g\in G$ for the latter.

The second challenge is analytic
and comes from the singularity of
the initial condition in $C^{-\frac12-\kappa}$ for $\kappa>0$ (this was already encountered in 2D in a more mild form, see Remark~\ref{rem:input_dist_2D}).
This singularity means, for example, that $\mcP A(0) \partial \Psi$ is ill-defined for generic
distributions $\Psi\in C^{-\frac12-\kappa}(\R\times\T^3 )$ and $A(0)\in C^{-\frac12-\kappa}(\T^3)$.
Similar to the discussion in Section~\ref{subsubsec:1st_half}, this type of product appears in the Picard iteration used to solve~\eqref{SYM}.
As in Remark~\ref{rem:input_dist_2D},
this issue is addressed
by decomposing $A = \CP A(0) + \Psi + \hat A$, where $\Psi$ solves the SHE $\partial_t\Psi = \Delta \Psi + \xi^\eps$,
and solving for the `remainder' $\hat A$.
One then shows with separate stochastic bounds that $\CP A(0)\partial \Psi$ and $\Psi\partial \CP A(0)$ converge in $C^{-2-\kappa}(\R\times \T^3)$ as $\eps\downarrow0$.

A closely related issue \textit{not present} in 2D is that of
restarting the equation at some positive time $\tau>0$ to obtain maximal solutions.
This is because, for $\eps>0$, $A(\tau)$ and $\xi^\eps\restr_{[\tau,\infty)}$ see each other on a time interval of order $\eps^2$.
Since the regularities of $A(\tau)$ and $\partial\Psi$ add up to $<-2$,
this breaks the argument used to show that $\CP A(0) \partial \Psi$ and $\Psi\partial\CP A(0)$ converge as $\eps\downarrow0$ when $A(0)$ is independent of $\xi$.
To restart the equation, one instead leverages that $A(\tau)$ for $\tau>0$ is not a generic element of $\state$ but takes the form $A(\tau) = \Psi(\tau) + R(\tau)$ where $R(\tau)\in C^{-\kappa}(\T^3)$.
Since $\Psi$ is defined globally in time,
this decomposition allows one to restart the equation using the `generalised Da Prato--Debussche trick' from~\cite{BCCH21}.

\subsection{Gauge covariance}
\label{subsec:3D_gauge_covar}

Finally, we describe the proof of Theorem~\ref{thm:generative} in 3D.
The proof of Theorem~\ref{thm:generative}\ref{pt:gen_exists}
is similar to its 2D counterpart.
The only appreciable difference is that the measurable selection $S\colon\mfO\to\state$
is replaced by a Borel map $S\colon \state\to\state$
which preserves gauge orbits and such that
$\bar\Sigma(S(X)) \leq 2\inf_{Y\sim X} \bar \Sigma(Y)$ whenever the right-hand side is finite.
Here $\bar \Sigma\geq \Sigma$ is defined analogously to $\Sigma$ but with a stronger set of parameters $(\eta,\beta,\delta,\alpha,\theta)$. 
This complication is due to a lack of any nice known properties of $\mfO=\state/{\sim}$ in 3D (e.g.\ Polishness), see Remarks~\ref{rem:quotient_top} and~\ref{rem:no_gauge_group}, and we instead leverage compactness of the embedding $(\bar\state,\bar\Sigma)\hookrightarrow(\state,\Sigma)$.

The proof of Theorem~\ref{thm:generative}\ref{pt:covar}
is where we start to see a difference with the 2D case.
First, admitting for now Theorem~\ref{thm:B_bar_A},
we aim to prove that solutions to~\eqref{SYM} with bare mass $\check C$ and gauge equivalent initial conditions can be suitably coupled.
Namely, one has the following result.

\begin{lemma}[Coupling]\label{lem:coupling}
Suppose $A$ solves~\eqref{SYM} with bare mass $\check C$ and initial condition $a$.
Then, for any $b\sim a$,
there exists on the same probability space
a white noise $\bar \xi$ and a process $(B, g)$ such that
$g(t)\in \mfG^{\frac32-\kappa}$ and $A(t)^{g(t)}=B(t)$ for all $t>0$ (before blow-up of $A,B$)
and such that $B$ solves~\eqref{SYM} driven by $\bar\xi$ with bare mass $\check C$.
\end{lemma}

If $a^{g(0)}=b$ for some $g(0)\in \mfG^{\rho}$ for $\rho>\frac12$ as in Section~\ref{subsubsec:2nd_half},
then this result follows almost immediately from Theorem~\ref{thm:B_bar_A}.
However, unlike the 2D case, it is now possible that $b\sim a$ but no $g(0)$ exists such that $a^{g(0)}=b$, which leads to trouble in applying
Theorem~\ref{thm:B_bar_A} -- we effectively have no initial condition for $g$ in the PDE~\eqref{eq:SPDE_B}.

To circumvent the issue, we regularise $a$ and $b$ using the YM flow $\mcF_{t}$,
so that $\mcF_{t}(a)^{g_t} = \mcF_t(b)$  for all $t$ sufficiently small and some smooth $g_t$.
Then
$\mcF_{t}(a) \to a$ and $\mcF_{t}(b)\to b$ in $\state$
as $t\downarrow0$,
which in particular implies that $\limsup_{t\downarrow0}|g_t|_{C^\nu} < \infty$ for some $\nu\in(0,\frac12)$ due to Theorem~\ref{thm:g_bound}.
Therefore, there exists $g(0)\in \mfG^{\nu}$ such that $g_t\to g(0)$ in $\mfG^{\nu/2}$ along a subsequence.
We can then rewrite the equation for $g$ in~\eqref{eq:SPDE_B} in terms of $A=B^{g^{-1}}$, namely
\begin{equ}\label{eq:g_wrt_A}
g^{-1}\partial_t g = \partial_j(g^{-1}\partial_j g) + [A_j,g^{-1}\partial_jg]\;,
\end{equ}
and show that the system for $(A,g)$, where $A$ solves~\eqref{SYM}, is well-posed for any initial condition in $\state\times\mfG^\nu$ (vs. $\state\times \mfG^\rho$ with $\rho>\frac12$ for~\eqref{eq:SPDE_B} due to the multiplicative noise $\Ad_g \xi$).
One can then use continuity of $(A,g)$ with respect to initial conditions, the fact that $g$
takes values in $\mfG^{\frac32-\kappa}$ for positive times, and the joint continuity of the group action $\mfG^{\frac32-\kappa}\times\state \to \state$,
to prove Lemma~\ref{lem:coupling}.
The $g$ in Lemma~\ref{lem:coupling} solves precisely~\eqref{eq:g_wrt_A}
with initial condition $g(0)$.

With Lemma~\ref{lem:coupling} in hand, together with the fact that $g$ in its statement cannot blow up before $\Sigma(A,0)+\Sigma(B,0)$ blows up (again due to Theorem~\ref{thm:g_bound}),
it is relatively straightforward to prove Theorem~\ref{thm:generative}\ref{pt:covar} like we did in the 2D case.
See in particular the discussion around Figure~\ref{fig:Markov}.

It remains to prove Theorem~\ref{thm:B_bar_A} to complete the above argument.
Fix henceforth a non-anticipative mollifier $\moll$.
Like in Section~\ref{subsec:covar_2D},
it is natural consider the systems~\eqref{eq:SPDE_B_bares} and~\eqref{eq:SPDE_bar_A_bares},
and show that their renormalised versions converge to the same limit.
To do this, we show that the equations for $g$ and $\bar g$, which are now singular, do not require renormalisation,
and that the renormalised equations for $B$ and $\bar A$ are of exactly the same form as in Section~\ref{subsec:covar_2D}, namely
\begin{equ}\label{eq:B_3D_renorm}
\partial_t B = \Delta B + B\partial B + B^3 + \Ad_{g} \xi^\eps + (C^\eps_\BPHZ +\mathring C_1)B + (\tilde C^\eps + \mathring C_2) (\mrd g) g^{-1}\;,
\end{equ}
and
\begin{equ}\label{eq:bar_A_3D_renorm}
\partial_t \bar A = \Delta \bar A + \bar A\partial \bar A + \bar A^3 + \moll^\eps*(\Ad_{\bar g}\xi) + (C^\eps_\BPHZ +\mathring C_1)\bar A + (\tilde C^{0,\eps} + \mathring C_2) (\mrd \bar g )\bar g^{-1}\;.
\end{equ}
The constants $\tilde C^\eps, \tilde C^{0,\eps}$ are respectively called $C^\eps_{\Gauge}, C^{0,\eps}_{\Gauge}$ in~\cite{CCHS22},
and correspond to $\lambda \tilde C^\eps, \lambda \tilde C^{0,\eps}$ in Section~\ref{subsec:covar_2D};
$C^\eps_\BPHZ$ is the BPHZ constant in Theorem~\ref{thm:local} and corresponds to
$\lambda C^\eps_\SYM$ in Section~\ref{sec:2D}.

Both of these facts are again proven by introducing the new variables $h=(\mrd g) g^{-1}$ and $U=\Ad_g$ and writing the corresponding equations for $h,U$.
This time, instead of a direction computation with just two trees as in~\eqref{eq:2D_finite_trees},
a more complicated strategy relying on power counting and symmetry arguments is necessary.
The main insight which helps with this argument is that 
the trees appearing in the `lifted' $B$ and $\bar A$ equations are obtained by attaching (or grafting) the trees
appearing in the $U$ equation onto the trees appearing in the~\eqref{SYM} equation (see~\cite[Sec.~6.2.1]{CCHS22}).

The analytic theory for the $B$ and $\bar A$ equations is somewhat more involved than what we saw in Sections~\ref{subsec:covar_2D} and~\ref{subsec:local_sol_3D},
but the general strategy is the same:
we decompose the solution into the initial condition, globally defined singular terms, and a better behaved remainder, and solve for the last of these.
Restarting these equations, specifically the equation for $\bar A$, is also not straightforward
because, unlike in Section~\ref{subsec:local_sol_3D},
we are outside the scope of the generalised Da Prato--Debussche trick of~\cite{BCCH21}
since the multiplicative noise means that the most singular part of the solution
is not just the SHE.
Therefore, an entirely separate fixed point problem needs to be written for the restarted equation (see~\cite[Sec.~6.6]{CCHS22})
which leverages that the new initial condition comes from a modelled distribution defined for earlier times.
The proof that $B$ and $\bar A$ converge to the same limit uses the same $\eps$-dependent norms on regularity structures introduced in~\cite{CCHS20}.

To conclude the proof of Theorem~\ref{thm:B_bar_A}, like in Section~\ref{subsec:covar_2D}, it suffices to show that
\begin{equ}\label{eq:check_C_3D}
\check C \eqdef \lim_{\eps\downarrow0}(\tilde C^\eps - C^\eps_\BPHZ - \tilde C^{0,\eps}) \;\;
\text{exists and is independent of
non-anticipative $\moll$}\;.
\end{equ}
The difficulty compared to 2D is that we have no explicit formulae for the constants on the right-hand side.
Instead, the argument in~\cite{CCHS22} proceeds in two steps.

First, one shows that
\begin{equ}\label{eq:BPHZ_diff}
\limsup_{\eps\downarrow0}|\tilde C^\eps - C^\eps_\BPHZ|<\infty
\;\;\text{ and }\;\;
\limsup_{\eps\downarrow0}|\tilde C^{0,\eps}| < \infty\;.
\end{equ}
Arguing by contradiction, suppose that
\begin{equ}\label{eq:contradiction1}
|\tilde C^\eps - C^\eps_\BPHZ| \to\infty
\end{equ}
along a subsequence of $\eps\downarrow0$.
The idea, inspired by~\cite{BGHZ21}, is to introduce a parameter $\sigma^\eps\in \R$ in front of the noise which we take to zero as $\eps\downarrow0$ in a precise way.

Letting $\tilde C^\eps_{\sigma}, C^\eps_{\BPHZ,\sigma}$, etc. denote the 
renormalisation constants associated to the noise $\sigma\xi$,
it is not difficult to see that these constants depend polynomially on $\sigma$ and converge to zero as $\sigma\downarrow0$.
Therefore, by assumption~\eqref{eq:contradiction1},
we can send $\sigma^\eps \downarrow 0$ as $\eps\downarrow0$ in such a way that, after passing to another subsequence, we can find $\hat C\neq 0$ such that
\begin{equ}
\hat C^\eps \eqdef C^\eps_{\BPHZ,\sigma^\eps}-\tilde C^\eps_{\sigma^\eps} \to \hat C\;.
\end{equ}
Now we consider the equation for $B$ with bare masses $\mathring C_1=0$ and $\mathring C_2=\hat C^\eps$,
that is,
\begin{equ}
\partial_t B = \Delta B + B \partial B + B^3 + C^\eps_{\BPHZ,\sigma^\eps} B + \sigma^\eps\Ad_g\xi^\eps + (\tilde C^\eps_{\sigma^\eps} + \hat C^\eps)(\mrd g) g^{-1}\;.
\end{equ}
Remark that $\tilde C^\eps_{\sigma^\eps} + \hat C^\eps = C^\eps_{\BPHZ,\sigma^\eps}$ by definition,
and thus $B=A^g$ pathwise
where $A$ solves~\eqref{SYM} with bare mass zero, i.e.
\begin{equ}
\partial_t A = \Delta A + A \partial A + A^3 + C^\eps_{\BPHZ,\sigma^\eps} A + \sigma^\eps\xi^\eps\;.
\end{equ}
Since $\sigma^\eps\downarrow0$, $A$ converges as $\eps\downarrow0$ to solution of the deterministic YM flow (with DeTurck term)
\begin{equ}\label{eq:YM_flow_lim}
\partial_t A = \Delta A + A \partial A + A^3
\end{equ}
with the same initial condition.
Furthermore, because $\hat C^\eps \to \hat C$,
we can use joint continuity in the bare mass and the noise to see that
$B$ converges as $\eps\downarrow0$
to
the solution of
\begin{equ}\label{eq:B_eq_1}
\partial_t B = \Delta B + B \partial B + B^3 +  \hat C(\mrd g) g^{-1}\;.
\end{equ}
On the other hand, the equality $A^g=B$ is preserved under the $\eps\downarrow0$ limit, where $g$ solves~\eqref{eq:SPDE_B},
but this implies by the calculation in Section~\ref{subsec:gauge_covar} with $\xi\equiv 0$ that $B$ should also solve
\begin{equ}\label{eq:B_eq_2}
\partial_t B = \Delta B + B \partial B + B^3\;.
\end{equ}
This is clearly a contradiction because, for any $\hat C\neq0$, we can find initial conditions $(B(0),g(0))$
for which the solutions to~\eqref{eq:B_eq_1} and~\eqref{eq:B_eq_2} are different.

The argument that $\limsup_{\eps\downarrow0}|\tilde C^{0,\eps}| < \infty$ is similar.
Namely, if $|\tilde C^{0,\eps}| \to\infty$ along a subsequence,
then we can find $\sigma^\eps\downarrow0$ such that $\tilde C^{0,\eps}_{\sigma^\eps} \to \hat C \neq 0$ along another subsequence.
Now we consider
\begin{equ}
\partial \bar A
= \Delta \bar A + \bar A\partial\bar A + \bar A^3 + (C^\eps_{\BPHZ,\sigma^\eps} + \mathring C_1)\bar A
+ (\tilde C^{0,\eps}_{\sigma^\eps}+\mathring C_2 )(\mrd \bar g) \bar g^{-1} + \moll^\eps*(\Ad_{\bar g} \sigma^\eps\xi)
\end{equ}
with bare masses $\mathring C_1=0$ and $\mathring C_2=-\tilde C^{0,\eps}_{\sigma^\eps}$.
With this choice, the term $(\mrd \bar g) \bar g^{-1}$ vanishes and the solution is, by It{\^o} isometry (which requires $\moll$ non-anticipative),
equal in law to the solution of
\begin{equ}
\partial \tilde A = \Delta \tilde A + \tilde A\partial \tilde A + \tilde A^3 + C^\eps_{\BPHZ,\sigma^\eps} \tilde A + \sigma^\eps\xi^\eps\;.
\end{equ}
But now, since $\sigma^\eps\downarrow0$, $\tilde A$ converges in law again to the solution of the deterministic YM flow~\eqref{eq:YM_flow_lim},
while $\bar A$ converges to the solution of
\begin{equ}
\partial \bar A = \Delta \bar A + \bar A\partial\bar A + \bar A^3 - \hat C (\mrd \bar g) \bar g^{-1} \;,
\end{equ}
which are not equal for generic initial conditions,
bringing us to a contradiction.

\begin{remark}
The bounds~\eqref{eq:BPHZ_diff}
are actually used in the short time analysis of~\eqref{eq:B_3D_renorm}-\eqref{eq:bar_A_3D_renorm} since they allow us to relate
$B$ and $\bar A$ to a simpler equation with \textit{additive} noise through the above mechanism ($B$ is related through pathwise gauge transformations, $\bar A$ is related through equality in law, see~\cite[Sec.~6.6]{CCHS22}).
\end{remark}

To finish the proof of~\eqref{eq:check_C_3D}, it is left to show that
\begin{equ}\label{eq:lim_exists}
\lim_{\eps\downarrow0}  C^\eps_\BPHZ-\tilde C^\eps \;\;\text{exists and is independent of non-anticipative $\moll$}\;,
\end{equ}
and similarly for $\tilde C^{0,\eps}$.
Arguing again by contradiction,
suppose we have two subsequences $\eps\downarrow0$
along which $C^\eps_\BPHZ-\tilde C^\eps$ converges to distinct limits $\hat C_1$ and $\hat C_2$.
We then consider the $B$ equation with bare masses $\mathring C_1=0$ and $\mathring C_2 = C^\eps_\BPHZ - \tilde C^\eps$.
With this choice,
$B=A^g$ for all $\eps>0$ where $A$ solves~\eqref{SYM} with zero bare mass and $g$ solves~\eqref{eq:g_wrt_A}.
Since the limit of $(A,g)$ is the same along the two subsequences,
it follows that the limit of $B$ along the two subsequence is also the same almost surely.
However, it is possible to show that
the limiting bare mass constants $\hat C_1$ and $\hat C_2$ can be recovered from the $\eps\downarrow0$ limit of $(B,g)$, see~\cite[Appendix~D]{CCHS22},
which leads to a contradiction since we assumed $\hat C_1\neq \hat C_2$.
The same argument shows that the limit in~\eqref{eq:lim_exists} is independent of $\moll$ because the solution $(A,g)$ to~\eqref{SYM} and~\eqref{eq:g_wrt_A} is independent of $\moll$.

Finally, to prove that $\lim_{\eps\downarrow0} \tilde C^{0,\eps}$ exists and is independent of the non-anticipative mollifier $\moll$,
one argues in a similar way except, as earlier, we appeal to equality with~\eqref{SYM} \textit{in law}
and use instead that two limiting solutions $(\bar A,\bar g)$ with different bare masses cannot be equal \textit{in law}.

\begin{remark}
These final statements mimic
exactly what we saw in Section~\ref{subsec:covar_2D} where $\tilde C^{0,\eps}=0$ for non-anticipative $\moll$ and $\lim_{\eps\downarrow0} (\tilde C^\eps - C^\eps_\BPHZ)$ is given by~\eqref{eq:check_C_2D}.
By analogy,~\eqref{eq:lim_exists} should hold for any mollifier, not necessarily non-anticipative,
but this is not necessarily the case for $\lim_{\eps\downarrow0}\tilde C^{0,\eps}$.
Furthermore, as in 2D, we expect that $\lim_{\eps\downarrow0}\tilde C^{0,\eps}=0$ for non-anticipative $\moll$.
\end{remark}

\section{Open problems}

We close with several open problems which we believe to be of interest.

\begin{itemize}
\item Does the Markov process on gauge orbits in Theorem~\ref{thm:generative}
possess a unique invariant measure?
Existence of the invariant measure should imply uniqueness due to 
the strong Feller property~\cite{HM18_Feller} and full support theorem for SPDEs~\cite{HS22_Support}.
Furthermore, the invariant measure for $d=2$ is expected to be the YM measure associated to the trivial principal $G$-bundle on $\T^2$ constructed in~\cite{Sengupta97,Levy03,Levy06}.
For $d=3$, this would provide the first construction of the YM measure in 3D, even in finite volume.

\item Can the analysis in~\cite{CCHS20,CCHS22} be extended to infinite volume $\T^d\rightsquigarrow \R^d$?
This is non-trivial even for $d=2$, although the YM measure on $\R^2$ is arguably simpler~\cite{Driver89}.

\item Can one extend these results beyond the case that the underlying principal bundle $P\to \T^d$ is trivial?
For non-trivial principal bundles, one can no longer write connections as globally defined $1$-forms, which complicates the solution theory.

\item For $d=3$, can one modify the construction of the state 
space $\state$ in Section~\ref{subsec:3D_state}
so that the gauge equivalence $\sim$ is determined by a gauge group or a similar structure?
This would yield a notion of gauge equivalence conceptually closer to the classical space of gauge orbits
and would carry a number of technical advantages (see Remarks~\ref{rem:quotient_top} and~\ref{rem:no_gauge_group} and the start of Section~\ref{subsec:3D_gauge_covar}).

\item Taking $G$ to be one of the classical groups, e.g.\ $G=U(N)$, what is the behaviour of the dynamic as $N\to\infty$?
In 2D, the associated YM measure is known to converge to a deterministic object called the master field~\cite{Levy17,DN20,Dahlqvist_Lemoine_22_I,Dahlqvist_Lemoine_22_II}, which is governed by the Makeenko--Migdal equations~\cite{MM79};
see~\cite{Levy20} for a survey.
No such result is rigorously known in 3D (the measure at finite $N$ has not been constructed).
It would be interesting if one can use stochastic quantisation to recover some of the known results in 2D and obtain new results in 3D;
see~\cite{SSZ22} where the Langevin dynamic is used to derive the finite $N$ master loop equation on the lattice.
\end{itemize}



\subsubsection*{Acknowledgements}

It is a great pleasure to thank Ajay Chandra, Martin Hairer, and Hao Shen for many interesting discussions and explorations while carrying out the reviewed works.

\bibliographystyle{./Martin}
\bibliography{./refs}

\newcommand{\etalchar}[1]{$^{#1}$}
\def\cprime{$'$} \def\polhk#1{\setbox0=\hbox{#1}{\ooalign{\hidewidth
  \lower1.5ex\hbox{`}\hidewidth\crcr\unhbox0}}}
\begin{thebibliography}{CCHS22b}
\def\myhref#1#2{\href{#2}{\nolinkurl{#1}}}

\bibitem[AK20]{AK20}
\textsc{S.~Albeverio} and \textsc{S.~Kusuoka}.
\newblock The invariant measure and the flow associated to the
  {$\Phi^4_3$}-quantum field model.
\newblock \emph{Ann. Sc. Norm. Super. Pisa Cl. Sci. (5)} \textbf{20}, no.~4,
  (2020), 1359--1427.
\newblock
  \myhref{doi:10.2422/2036-2145.201809_008}{https://dx.doi.org/10.2422/2036-2145.201809_008}.

\bibitem[Bal85]{Balaban85IV}
\textsc{T.~Balaban}.
\newblock Ultraviolet stability of three-dimensional lattice pure gauge field
  theories.
\newblock \emph{Comm. Math. Phys.} \textbf{102}, no.~2, (1985), 255--275.
\newblock
  \myhref{doi:10.1007/BF01229380}{https://dx.doi.org/10.1007/BF01229380}.

\bibitem[Bal87]{Balaban87}
\textsc{T.~Balaban}.
\newblock Renormalization group approach to lattice gauge field theories. {I}.
  {G}eneration of effective actions in a small field approximation and a
  coupling constant renormalization in four dimensions.
\newblock \emph{Comm. Math. Phys.} \textbf{109}, no.~2, (1987), 249--301.
\newblock
  \myhref{doi:10.1007/BF01215223}{https://dx.doi.org/10.1007/BF01215223}.

\bibitem[Bal89]{Balaban89}
\textsc{T.~Balaban}.
\newblock Large field renormalization. {II}. {L}ocalization, exponentiation,
  and bounds for the {R} operation.
\newblock \emph{Comm. Math. Phys.} \textbf{122}, no.~3, (1989), 355--392.
\newblock
  \myhref{doi:10.1007/BF01238433}{https://dx.doi.org/10.1007/BF01238433}.

\bibitem[BCCH21]{BCCH21}
\textsc{Y.~Bruned}, \textsc{A.~Chandra}, \textsc{I.~Chevyrev}, and
  \textsc{M.~Hairer}.
\newblock Renormalising {SPDE}s in regularity structures.
\newblock \emph{J. Eur. Math. Soc. (JEMS)} \textbf{23}, no.~3, (2021),
  869--947.
\newblock \myhref{doi:10.4171/jems/1025}{https://dx.doi.org/10.4171/jems/1025}.

\bibitem[BGHZ21]{BGHZ21}
\textsc{Y.~Bruned}, \textsc{F.~Gabriel}, \textsc{M.~Hairer}, and
  \textsc{L.~Zambotti}.
\newblock Geometric stochastic heat equations.
\newblock \emph{J. Amer. Math. Soc.} \textbf{35}, no.~1, (2021), 1--80.
\newblock \myhref{doi:10.1090/jams/977}{https://dx.doi.org/10.1090/jams/977}.

\bibitem[BHST87]{BHST87II}
\textsc{Z.~Bern}, \textsc{M.~B. Halpern}, \textsc{L.~Sadun}, and
  \textsc{C.~Taubes}.
\newblock Continuum regularization of quantum field theory. {II}. {G}auge
  theory.
\newblock \emph{Nuclear Phys. B} \textbf{284}, no.~1, (1987), 35--91.
\newblock
  \myhref{doi:10.1016/0550-3213(87)90026-5}{https://dx.doi.org/10.1016/0550-3213(87)90026-5}.

\bibitem[BHZ19]{BHZ16}
\textsc{Y.~Bruned}, \textsc{M.~Hairer}, and \textsc{L.~Zambotti}.
\newblock Algebraic renormalisation of regularity structures.
\newblock \emph{Invent. Math.} \textbf{215}, no.~3, (2019), 1039--1156.
\newblock
  \myhref{doi:10.1007/s00222-018-0841-x}{https://dx.doi.org/10.1007/s00222-018-0841-x}.

\bibitem[CC21a]{Sourav_state}
\textsc{S.~{Cao}} and \textsc{S.~{Chatterjee}}.
\newblock {A state space for 3D Euclidean Yang-Mills theories}.
\newblock \emph{arXiv e-prints} (2021).
\newblock \myhref{arXiv:2111.12813}{https://arxiv.org/abs/2111.12813}.

\bibitem[CC21b]{Sourav_flow}
\textsc{S.~{Cao}} and \textsc{S.~{Chatterjee}}.
\newblock {The Yang-Mills heat flow with random distributional initial data}.
\newblock \emph{arXiv e-prints} (2021).
\newblock \myhref{arXiv:2111.10652}{https://arxiv.org/abs/2111.10652}.

\bibitem[CCHS22a]{CCHS20}
\textsc{A.~{Chandra}}, \textsc{I.~{Chevyrev}}, \textsc{M.~{Hairer}}, and
  \textsc{H.~{Shen}}.
\newblock {Langevin dynamic for the 2D Yang-Mills measure}.
\newblock \emph{Publ. Math. Inst. Hautes \'{E}tudes Sci.} (2022).
\newblock \myhref{arXiv:2006.04987}{https://arxiv.org/abs/2006.04987}.
\newblock
  \myhref{doi:10.1007/s10240-022-00132-0}{https://dx.doi.org/10.1007/s10240-022-00132-0}.

\bibitem[CCHS22b]{CCHS22}
\textsc{A.~{Chandra}}, \textsc{I.~{Chevyrev}}, \textsc{M.~{Hairer}}, and
  \textsc{H.~{Shen}}.
\newblock {Stochastic quantisation of Yang-Mills-Higgs in 3D}.
\newblock \emph{arXiv e-prints} (2022).
\newblock \myhref{arXiv:2201.03487}{https://arxiv.org/abs/2201.03487}.

\bibitem[CG13]{CG13}
\textsc{N.~Charalambous} and \textsc{L.~Gross}.
\newblock The {Y}ang-{M}ills heat semigroup on three-manifolds with boundary.
\newblock \emph{Comm. Math. Phys.} \textbf{317}, no.~3, (2013), 727--785.
\newblock
  \myhref{doi:10.1007/s00220-012-1558-0}{https://dx.doi.org/10.1007/s00220-012-1558-0}.

\bibitem[CH16]{CH16}
\textsc{A.~{Chandra}} and \textsc{M.~{Hairer}}.
\newblock {An analytic BPHZ theorem for regularity structures}.
\newblock \emph{ArXiv e-prints} (2016).
\newblock \myhref{arXiv:1612.08138}{https://arxiv.org/abs/1612.08138}.

\bibitem[Cha19]{Chatterjee18}
\textsc{S.~Chatterjee}.
\newblock Yang-{M}ills for probabilists.
\newblock In \emph{Probability and analysis in interacting physical systems},
  vol. 283 of \emph{Springer Proc. Math. Stat.},  1--16. Springer, Cham, 2019.
\newblock
  \myhref{doi:10.1007/978-3-030-15338-0_1}{https://dx.doi.org/10.1007/978-3-030-15338-0_1}.

\bibitem[Che19]{Chevyrev19YM}
\textsc{I.~Chevyrev}.
\newblock Yang-{M}ills measure on the two-dimensional torus as a random
  distribution.
\newblock \emph{Comm. Math. Phys.} \textbf{372}, no.~3, (2019), 1027--1058.
\newblock
  \myhref{doi:10.1007/s00220-019-03567-5}{https://dx.doi.org/10.1007/s00220-019-03567-5}.

\bibitem[{Che}22]{Chevyrev22_norm_inf}
\textsc{I.~{Chevyrev}}.
\newblock {Norm inflation for a non-linear heat equation with Gaussian initial
  conditions}.
\newblock \emph{arXiv e-prints} (2022).
\newblock \myhref{arXiv:2205.14350}{https://arxiv.org/abs/2205.14350}.

\bibitem[CW17]{ChandraWeber17}
\textsc{A.~Chandra} and \textsc{H.~Weber}.
\newblock Stochastic {PDE}s, regularity structures, and interacting particle
  systems.
\newblock \emph{Ann. Fac. Sci. Toulouse Math. (6)} \textbf{26}, no.~4, (2017),
  847--909.
\newblock \myhref{doi:10.5802/afst.1555}{https://dx.doi.org/10.5802/afst.1555}.

\bibitem[DDPR13]{MR3133916}
\textsc{L.~Del~Debbio}, \textsc{A.~Patella}, and \textsc{A.~Rago}.
\newblock Space-time symmetries and the {Y}ang-{M}ills gradient flow.
\newblock \emph{J. High Energy Phys.} \textbf{2013}, no. 212(2013).
\newblock
  \myhref{doi:10.1007/JHEP11(2013)212}{https://dx.doi.org/10.1007/JHEP11(2013)212}.

\bibitem[DeT83]{deturck83}
\textsc{D.~M. DeTurck}.
\newblock Deforming metrics in the direction of their {R}icci tensors.
\newblock \emph{J. Differential Geom.} \textbf{18}, no.~1, (1983), 157--162.
\newblock
  \myhref{doi:10.4310/JDG/1214509286}{https://dx.doi.org/10.4310/JDG/1214509286}.

\bibitem[DL22a]{Dahlqvist_Lemoine_22_II}
\textsc{A.~{Dahlqvist}} and \textsc{T.~{Lemoine}}.
\newblock {Large N limit of the Yang-Mills measure on compact surfaces II:
  Makeenko-Migdal equations and planar master field}.
\newblock \emph{arXiv e-prints} (2022).
\newblock \myhref{arXiv:2201.05886}{https://arxiv.org/abs/2201.05886}.

\bibitem[DL22b]{Dahlqvist_Lemoine_22_I}
\textsc{A.~{Dahlqvist}} and \textsc{T.~{Lemoine}}.
\newblock {Large N limit of Yang-Mills partition function and Wilson loops on
  compact surfaces}.
\newblock \emph{arXiv e-prints} (2022).
\newblock \myhref{arXiv:2201.05882}{https://arxiv.org/abs/2201.05882}.

\bibitem[DN20]{DN20}
\textsc{A.~Dahlqvist} and \textsc{J.~R. Norris}.
\newblock Yang-{M}ills measure and the master field on the sphere.
\newblock \emph{Comm. Math. Phys.} \textbf{377}, no.~2, (2020), 1163--1226.
\newblock
  \myhref{doi:10.1007/s00220-020-03773-6}{https://dx.doi.org/10.1007/s00220-020-03773-6}.

\bibitem[Don85]{Donaldson}
\textsc{S.~K. Donaldson}.
\newblock Anti self-dual {Y}ang-{M}ills connections over complex algebraic
  surfaces and stable vector bundles.
\newblock \emph{Proc. London Math. Soc. (3)} \textbf{50}, no.~1, (1985), 1--26.
\newblock
  \myhref{doi:10.1112/plms/s3-50.1.1}{https://dx.doi.org/10.1112/plms/s3-50.1.1}.

\bibitem[Dri87]{Driver87}
\textsc{B.~K. Driver}.
\newblock Convergence of the {${\rm U}(1)_4$} lattice gauge theory to its
  continuum limit.
\newblock \emph{Comm. Math. Phys.} \textbf{110}, no.~3, (1987), 479--501.
\newblock
  \myhref{doi:10.1007/BF01212424}{https://dx.doi.org/10.1007/BF01212424}.

\bibitem[Dri89]{Driver89}
\textsc{B.~K. Driver}.
\newblock Y{M{${}_2$}}: continuum expectations, lattice convergence, and
  lassos.
\newblock \emph{Comm. Math. Phys.} \textbf{123}, no.~4, (1989), 575--616.
\newblock
  \myhref{doi:10.1007/BF01218586}{https://dx.doi.org/10.1007/BF01218586}.

\bibitem[{Duc}21]{Duch21}
\textsc{P.~{Duch}}.
\newblock {Flow equation approach to singular stochastic PDEs}.
\newblock \emph{arXiv e-prints} (2021).
\newblock \myhref{arXiv:2109.11380}{https://arxiv.org/abs/2109.11380}.

\bibitem[Fed86]{Federbush86}
\textsc{P.~Federbush}.
\newblock A phase cell approach to {Y}ang-{M}ills theory. {I}. {M}odes,
  lattice-continuum duality.
\newblock \emph{Comm. Math. Phys.} \textbf{107}, no.~2, (1986), 319--329.
\newblock
  \myhref{doi:10.1007/BF01209397}{https://dx.doi.org/10.1007/BF01209397}.

\bibitem[FH20]{FrizHairer20}
\textsc{P.~K. Friz} and \textsc{M.~Hairer}.
\newblock \emph{A course on rough paths}.
\newblock Universitext. Springer, Cham, 2020,  xvi+346.
\newblock With an introduction to regularity structures, Second edition of [
  3289027].
\newblock
  \myhref{doi:10.1007/978-3-030-41556-3}{https://dx.doi.org/10.1007/978-3-030-41556-3}.

\bibitem[FHK{\etalchar{+}}12]{Fodor12}
\textsc{Z.~Fodor}, \textsc{K.~Holland}, \textsc{J.~Kuti}, \textsc{D.~Nogradi},
  and \textsc{C.~H. Wong}.
\newblock The {Y}ang-{M}ills gradient flow in finite volume.
\newblock \emph{J. High Energy Phys.} \textbf{2012}, no.~11, (2012), 007, front
  matter + 16.
\newblock
  \myhref{doi:10.1007/JHEP11(2012)007}{https://dx.doi.org/10.1007/JHEP11(2012)007}.

\bibitem[Fin91]{Fine91}
\textsc{D.~S. Fine}.
\newblock Quantum {Y}ang-{M}ills on a {R}iemann surface.
\newblock \emph{Comm. Math. Phys.} \textbf{140}, no.~2, (1991), 321--338.
\newblock
  \myhref{doi:10.1007/BF02099502}{https://dx.doi.org/10.1007/BF02099502}.

\bibitem[GH21]{GH21}
\textsc{M.~Gubinelli} and \textsc{M.~Hofmanov\'{a}}.
\newblock A {PDE} construction of the {E}uclidean {$\phi_3^4$} quantum field
  theory.
\newblock \emph{Comm. Math. Phys.} \textbf{384}, no.~1, (2021), 1--75.
\newblock
  \myhref{doi:10.1007/s00220-021-04022-0}{https://dx.doi.org/10.1007/s00220-021-04022-0}.

\bibitem[GIP15]{GIP15}
\textsc{M.~Gubinelli}, \textsc{P.~Imkeller}, and \textsc{N.~Perkowski}.
\newblock Paracontrolled distributions and singular {PDE}s.
\newblock \emph{Forum Math. Pi} \textbf{3}, (2015), e6, 75.
\newblock \myhref{arXiv:1210.2684v3}{https://arxiv.org/abs/1210.2684v3}.
\newblock
  \myhref{doi:10.1017/fmp.2015.2}{https://dx.doi.org/10.1017/fmp.2015.2}.

\bibitem[GKS89]{GKS89}
\textsc{L.~Gross}, \textsc{C.~King}, and \textsc{A.~Sengupta}.
\newblock Two-dimensional {Y}ang-{M}ills theory via stochastic differential
  equations.
\newblock \emph{Ann. Physics} \textbf{194}, no.~1, (1989), 65--112.
\newblock
  \myhref{doi:10.1016/0003-4916(89)90032-8}{https://dx.doi.org/10.1016/0003-4916(89)90032-8}.

\bibitem[Gro83]{Gross83}
\textsc{L.~Gross}.
\newblock Convergence of {${\rm U}(1)\sb{3}$} lattice gauge theory to its
  continuum limit.
\newblock \emph{Comm. Math. Phys.} \textbf{92}, no.~2, (1983), 137--162.
\newblock
  \myhref{doi:10.1007/BF01210842}{https://dx.doi.org/10.1007/BF01210842}.

\bibitem[Hai14]{Hairer14}
\textsc{M.~Hairer}.
\newblock A theory of regularity structures.
\newblock \emph{Invent. Math.} \textbf{198}, no.~2, (2014), 269--504.
\newblock \myhref{arXiv:1303.5113}{https://arxiv.org/abs/1303.5113}.
\newblock
  \myhref{doi:10.1007/s00222-014-0505-4}{https://dx.doi.org/10.1007/s00222-014-0505-4}.

\bibitem[Hai16]{Hairer16_CDM}
\textsc{M.~Hairer}.
\newblock Regularity structures and the dynamical {$\Phi^4_3$} model.
\newblock In \emph{Current developments in mathematics 2014},  1--49. Int.
  Press, Somerville, MA, 2016.
\newblock \myhref{arXiv:1508.05261}{https://arxiv.org/abs/1508.05261}.

\bibitem[HM18]{HM18_Feller}
\textsc{M.~Hairer} and \textsc{J.~Mattingly}.
\newblock The strong {F}eller property for singular stochastic {PDE}s.
\newblock \emph{Ann. Inst. Henri Poincar\'{e} Probab. Stat.} \textbf{54},
  no.~3, (2018), 1314--1340.
\newblock
  \myhref{doi:10.1214/17-AIHP840}{https://dx.doi.org/10.1214/17-AIHP840}.

\bibitem[HS22a]{HS22_Support}
\textsc{M.~Hairer} and \textsc{P.~Sch{\"o}nbauer}.
\newblock The support of singular stochastic partial differential equations.
\newblock \emph{Forum of Mathematics, Pi} \textbf{10}, (2022), e1.
\newblock
  \myhref{doi:10.1017/fmp.2021.18}{https://dx.doi.org/10.1017/fmp.2021.18}.

\bibitem[HS22b]{HS22}
\textsc{M.~Hairer} and \textsc{R.~Steele}.
\newblock The {$\Phi_3^4$} {M}easure {H}as {S}ub-{G}aussian {T}ails.
\newblock \emph{J. Stat. Phys.} \textbf{186}, no.~3, (2022), Paper No. 38.
\newblock
  \myhref{doi:10.1007/s10955-021-02866-3}{https://dx.doi.org/10.1007/s10955-021-02866-3}.

\bibitem[HT04]{HT04}
\textsc{M.-C. Hong} and \textsc{G.~Tian}.
\newblock Global existence of the {$m$}-equivariant {Y}ang-{M}ills flow in four
  dimensional spaces.
\newblock \emph{Comm. Anal. Geom.} \textbf{12}, no. 1-2, (2004), 183--211.
\newblock
  \myhref{doi:10.4310/CAG.2004.v12.n1.a10}{https://dx.doi.org/10.4310/CAG.2004.v12.n1.a10}.

\bibitem[JW06]{JW06}
\textsc{A.~Jaffe} and \textsc{E.~Witten}.
\newblock Quantum {Y}ang-{M}ills theory.
\newblock In \emph{The millennium prize problems},  129--152. Clay Math. Inst.,
  Cambridge, MA, 2006.

\bibitem[Kup16]{Kupiainen16}
\textsc{A.~Kupiainen}.
\newblock Renormalization group and stochastic {PDE}s.
\newblock \emph{Ann. Henri Poincar\'{e}} \textbf{17}, no.~3, (2016), 497--535.
\newblock
  \myhref{doi:10.1007/s00023-015-0408-y}{https://dx.doi.org/10.1007/s00023-015-0408-y}.

\bibitem[L{\'e}v03]{Levy03}
\textsc{T.~L{\'e}vy}.
\newblock Yang-{M}ills measure on compact surfaces.
\newblock \emph{Mem. Amer. Math. Soc.} \textbf{166}, no. 790, (2003), xiv+122.
\newblock \myhref{arXiv:math/0101239}{https://arxiv.org/abs/math/0101239}.
\newblock \myhref{doi:10.1090/memo/0790}{https://dx.doi.org/10.1090/memo/0790}.

\bibitem[L{\'e}v06]{Levy06}
\textsc{T.~L{\'e}vy}.
\newblock Discrete and continuous {Y}ang-{M}ills measure for non-trivial
  bundles over compact surfaces.
\newblock \emph{Probab. Theory Related Fields} \textbf{136}, no.~2, (2006),
  171--202.
\newblock
  \myhref{doi:10.1007/s00440-005-0478-8}{https://dx.doi.org/10.1007/s00440-005-0478-8}.

\bibitem[L{\'e}v17]{Levy17}
\textsc{T.~L{\'e}vy}.
\newblock The master field on the plane.
\newblock \emph{Ast\'{e}risque} , no. 388, (2017), ix+201.
\newblock \myhref{arXiv:1112.2452}{https://arxiv.org/abs/1112.2452}.

\bibitem[L{\'e}v20]{Levy20}
\textsc{T.~L{\'e}vy}.
\newblock Two-dimensional quantum {Y}ang-{M}ills theory and the
  {M}akeenko-{M}igdal equations.
\newblock In \emph{Frontiers in analysis and probability---in the spirit of the
  {S}trasbourg-{Z}\"{u}rich meetings},  275--325. Springer, Cham, 2020.
\newblock
  \myhref{doi:10.1007/978-3-030-56409-4\\_7}{https://dx.doi.org/10.1007/978-3-030-56409-4\%5C_7}.

\bibitem[LS17]{LevySengupta17}
\textsc{T.~L\'{e}vy} and \textsc{A.~Sengupta}.
\newblock Four chapters on low-dimensional gauge theories.
\newblock In \emph{Stochastic geometric mechanics}, vol. 202 of \emph{Springer
  Proc. Math. Stat.},  115--167. Springer, Cham, 2017.
\newblock
  \myhref{doi:10.1007/978-3-319-63453-1_7}{https://dx.doi.org/10.1007/978-3-319-63453-1_7}.

\bibitem[L{\"u}s10]{Luscher10}
\textsc{M.~L{\"u}scher}.
\newblock Properties and uses of the {W}ilson flow in lattice {QCD}.
\newblock \emph{J. High Energy Phys.} \textbf{2010}, no.~8, (2010), 071, 18.
\newblock
  \myhref{doi:10.1007/JHEP08(2010)071}{https://dx.doi.org/10.1007/JHEP08(2010)071}.

\bibitem[Lyo94]{Lyons94}
\textsc{T.~Lyons}.
\newblock Differential equations driven by rough signals. {I}. {A}n extension
  of an inequality of {L}. {C}. {Y}oung.
\newblock \emph{Math. Res. Lett.} \textbf{1}, no.~4, (1994), 451--464.
\newblock
  \myhref{doi:10.4310/MRL.1994.v1.n4.a5}{https://dx.doi.org/10.4310/MRL.1994.v1.n4.a5}.

\bibitem[Mig75]{Migdal75}
\textsc{A.~A. Migdal}.
\newblock {Recursion Equations in Gauge Theories}.
\newblock \emph{Sov. Phys. JETP} \textbf{42}, (1975), 413.

\bibitem[MM79]{MM79}
\textsc{Y.~M. Makeenko} and \textsc{A.~A. Migdal}.
\newblock {Exact Equation for the Loop Average in Multicolor QCD}.
\newblock \emph{Phys. Lett. B} \textbf{88}, (1979), 135.
\newblock [Erratum: Phys.Lett.B 89, 437 (1980)].
\newblock
  \myhref{doi:10.1016/0370-2693(79)90131-X}{https://dx.doi.org/10.1016/0370-2693(79)90131-X}.

\bibitem[MRS93]{MRS93}
\textsc{J.~Magnen}, \textsc{V.~Rivasseau}, and \textsc{R.~S\'{e}n\'{e}or}.
\newblock Construction of {${\rm YM}_4$} with an infrared cutoff.
\newblock \emph{Comm. Math. Phys.} \textbf{155}, no.~2, (1993), 325--383.
\newblock
  \myhref{doi:10.1007/BF02097397}{https://dx.doi.org/10.1007/BF02097397}.

\bibitem[MW20]{MoinatWeber20}
\textsc{A.~Moinat} and \textsc{H.~Weber}.
\newblock Space-time localisation for the dynamic {$\Phi^4_3$} model.
\newblock \emph{Comm. Pure Appl. Math.} \textbf{73}, no.~12, (2020),
  2519--2555.
\newblock \myhref{doi:10.1002/cpa.21925}{https://dx.doi.org/10.1002/cpa.21925}.

\bibitem[NN06]{NN06}
\textsc{R.~Narayanan} and \textsc{H.~Neuberger}.
\newblock Infinite {$N$} phase transitions in continuum {W}ilson loop
  operators.
\newblock \emph{J. High Energy Phys.} \textbf{2006}, no.~3, (2006), 064, 32.
\newblock
  \myhref{doi:10.1088/1126-6708/2006/03/064}{https://dx.doi.org/10.1088/1126-6708/2006/03/064}.

\bibitem[OW19]{OW19}
\textsc{F.~Otto} and \textsc{H.~Weber}.
\newblock Quasilinear {SPDE}s via rough paths.
\newblock \emph{Arch. Ration. Mech. Anal.} \textbf{232}, no.~2, (2019),
  873--950.
\newblock
  \myhref{doi:10.1007/s00205-018-01335-8}{https://dx.doi.org/10.1007/s00205-018-01335-8}.

\bibitem[PW81]{ParisiWu}
\textsc{G.~Parisi} and \textsc{Y.~S. Wu}.
\newblock Perturbation theory without gauge fixing.
\newblock \emph{Sci. Sinica} \textbf{24}, no.~4, (1981), 483--496.
\newblock
  \myhref{doi:10.1360/ya1981-24-4-483}{https://dx.doi.org/10.1360/ya1981-24-4-483}.

\bibitem[Rad92]{Rade92}
\textsc{J.~Rade}.
\newblock On the {Y}ang-{M}ills heat equation in two and three dimensions.
\newblock \emph{J. Reine Angew. Math.} \textbf{431}, (1992), 123--163.
\newblock
  \myhref{doi:10.1515/crll.1992.431.123}{https://dx.doi.org/10.1515/crll.1992.431.123}.

\bibitem[{Sad}87]{Sadun}
\textsc{L.~A. {Sadun}}.
\newblock \emph{{Continuum Regularized Yang-Mills Theory}}.
\newblock Ph.D. thesis, University of California, Berkeley, CA (USA), 1987.

\bibitem[Sen97]{Sengupta97}
\textsc{A.~Sengupta}.
\newblock Gauge theory on compact surfaces.
\newblock \emph{Mem. Amer. Math. Soc.} \textbf{126}, no. 600, (1997), viii+85.
\newblock \myhref{doi:10.1090/memo/0600}{https://dx.doi.org/10.1090/memo/0600}.

\bibitem[She21]{Shen21}
\textsc{H.~Shen}.
\newblock Stochastic quantization of an {A}belian gauge theory.
\newblock \emph{Comm. Math. Phys.} \textbf{384}, no.~3, (2021), 1445--1512.
\newblock
  \myhref{doi:10.1007/s00220-021-04114-x}{https://dx.doi.org/10.1007/s00220-021-04114-x}.

\bibitem[SSZ22]{SSZ22}
\textsc{H.~{Shen}}, \textsc{S.~A. {Smith}}, and \textsc{R.~{Zhu}}.
\newblock {A new derivation of the finite $N$ master loop equation for lattice
  Yang-Mills}.
\newblock \emph{arXiv e-prints} (2022).
\newblock \myhref{arXiv:2202.00880}{https://arxiv.org/abs/2202.00880}.

\bibitem[You36]{Young}
\textsc{L.~C. Young}.
\newblock An inequality of the {H}\"older type, connected with {S}tieltjes
  integration.
\newblock \emph{Acta Math.} \textbf{67}, no.~1, (1936), 251--282.
\newblock
  \myhref{doi:10.1007/BF02401743}{https://dx.doi.org/10.1007/BF02401743}.

\bibitem[Zwa81]{zwanziger81}
\textsc{D.~Zwanziger}.
\newblock Covariant quantization of gauge fields without {G}ribov ambiguity.
\newblock \emph{Nuclear Phys. B} \textbf{192}, no.~1, (1981), 259--269.
\newblock
  \myhref{doi:10.1016/0550-3213(81)90202-9}{https://dx.doi.org/10.1016/0550-3213(81)90202-9}.

\end{thebibliography}

\end{document}